\documentclass[a4paper,twoside,11pt,leqno]{article}
\usepackage{amsthm}
\usepackage{amsmath}
\usepackage{natbib}
\RequirePackage[colorlinks]{hyperref}
\usepackage{amsfonts}

\def\authors{George Lowther}
\def\runauthor{\authors}
\def\head{A Generalized Backward Equation For One Dimensional Processes}
\def\runhead{A Generalized Backward Equation}
\def\keywords{diffusion, strong Markov, martingale, quasimartingale, Dirichlet process.}
\def\classifications{60J60, 60J25, 60G44, 60H10.}

\makeatletter
\def\@evenhead{\thepage\hfill{\small\MakeUppercase{\runauthor}}\hfill}
\def\@oddhead{\hfill{\small\MakeUppercase{\runhead}}\hfill\thepage}
\def\blfootnote{\xdef\@thefnmark{}\@footnotetext}

\makeatother

\newcommand{\halfplane}{\reals_+\times\reals}
\newcommand{\reals}{\mathbb{R}} 
\newcommand{\E}[1]{\mathbb{E}\left[#1\right]} 
\newcommand{\setsF}{\mathcal{F}} 
\newcommand{\setsG}{\mathcal{G}}                
\newcommand{\nicefunc}{\mathcal{D}}
\newcommand{\nicefuncK}{\mathcal{D}_{\rm K}}
\newcommand{\pd}[1]{_{\!,_{#1}}}
\newcommand{\dlimint}[4]{\int_{#1}^{#2\!\!\!}\!\int_{#3}^{#4}} 
\newcommand{\dblint}{\int\!\!\!\int} 
\newcommand{\Dom}{\mathcal D}
\newcommand{\nat}{\mathbb{N}}

\newcommand{\Prob}[1]{\mathbb{P}\left(#1\right)} 

\newcommand{\slbrack}{{[\hspace*{-0.15em}[}}  
\newcommand{\srbrack}{{]\hspace*{-0.15em}]}}  
\newcommand{\PP}{\mathbb{P}}

\newcommand\VarX[1]{{{\rm Var}_X(#1)}}

\newcommand{\zcqv}{{z.c.q.v.}}
\newcommand{\dotp}{{\cdot}} 
\newcommand{\ucp}{\xrightarrow{\rm\ ucp\ }} 
\newcommand{\cadlag}{c\`adl\`ag}
\newtheorem{definition}{Definition}[section]
\newtheorem{theorem}[definition]{Theorem}
\newtheorem{lemma}[definition]{Lemma}
\newtheorem{corollary}[definition]{Corollary}
\numberwithin{equation}{section}

\begin{document}

\title{\head}
\blfootnote{\emph{Key Words}: \keywords}
\blfootnote{\emph{AMS 2000 Classification}: \classifications}
\author{\authors}
\date{}
\maketitle
\thispagestyle{empty}

\begin{abstract}
Suppose that a real valued process $X$ is given as a solution to a stochastic differential equation. Then, for any twice continuously differentiable function $f$, the backward Kolmogorov equation gives a condition for $f(t,X)$ to be a local martingale.

We generalize the backward equation in two main ways. First, it is extended to non-differentiable functions. Second, the process $X$ is not required to satisfy an SDE. Instead, it is only required to be a quasimartingale satisfying an integrability condition, and the martingale condition for $f(t,X)$ is then expressed in terms of the marginal distributions, drift measure and jumps of $X$.
The proof involves the stochastic calculus of Dirichlet processes and a time-reversal argument.

These results are then applied to show that a continuous and strong Markov martingale is uniquely determined by its marginal distributions.
\end{abstract}

\maketitle

\section{Introduction}
\label{sec:Intro}

Suppose that we have a real valued process $X$ satisfying a stochastic differential equation of the form
\begin{equation}\label{eqn:SDE}
dX_t = \sigma(t,X_t)\,dW_t + b(t,X_t)\,dt
\end{equation}
where $W$ is a Brownian motion.
The uniqueness of solutions to this stochastic differential equation in the case where $\sigma$ and $b$ are Lipschitz continuous is well known (\citep{Rogers} Theorem 7.2). In fact, as we are just considering the case of one dimensional diffusions, the Yamada-Watanabe theorem (\citep{Rogers} Theorem 40.1) says that the Lipschitz condition for $\sigma$ can be relaxed to a H\"older condition of order $1/2$.
When $\sigma>0$ then a result of Stroock and Varadhan in \citep{Stroock3} shows that weak uniqueness holds as long as $\sigma$ is continuous and $b$ is locally bounded.

An alternative approach is to consider the backward Kolmogorov equation (or Feynman-Kac formula). Given a twice continuously differentiable function $f(t,x)$, this says that $f(t,X_t)$ is a local martingale if
\begin{equation}\label{eqn:Feynman Kac}
\frac{\partial}{\partial t}f(t,x)+\frac{1}{2}\sigma^2(t,x)\frac{\partial^2}{\partial x^2}f(t,x)+b(t,x)\frac{\partial}{\partial x}f(t,x)=0,
\end{equation}
which is a simple consequence of It\^o's formula.
Suppose that $\sigma$ and $b$ are sufficiently well behaved so that for any $T>0$ and smooth bounded function $g(x)$, there exists a bounded solution to (\ref{eqn:Feynman Kac}) for $t\le T$ and satisfying the boundary condition $f(T,x)=g(x)$. Then, $f(t,X_t)$ would be a martingale and consequently,
\begin{equation}\label{eqn:f is exp g}
f(t,X_t)=\E{g(X_T)|\setsF_t}.
\end{equation}
This shows that $X$ is Markov and determines its transition probabilities. So, the joint distribution of $X$ is uniquely defined once $X_0$ is known. Conversely, given a solution $X$ to (\ref{eqn:SDE}) we could define $f$ by equation (\ref{eqn:f is exp g}) and attempt to show that it is twice continuously differentiable, so it would satisfy the backward equation. Unfortunately this is not possible in general, in which case uniqueness of SDE (\ref{eqn:SDE}) can fail.

The aim of this paper is prove a generalization of the backward equation. Firstly, we generalize to include non-differentiable functions $f$. Secondly, we shall extend it to a much more general class of processes, which do not have to satisfy any stochastic differential equation such as (\ref{eqn:SDE}). In fact, we shall extend it to all quasimartingales satisfying an integrability condition. In that case the diffusion coefficients $\sigma$ and $b$ might not exist, and equation (\ref{eqn:Feynman Kac}) is replaced by one involving the marginal distributions, drift measure and jump components of $X$.

This result is applied to the example of continuous and strong Markov martingales, showing that they are uniquely determined by their marginal distributions.
Throughout this paper, we work with a filtered probability space $(\Omega,\setsF,(\setsF_t)_{t\in\reals_+},\PP)$, and the strong Markov property is defined as follows.
\begin{definition}\label{defn:strong markov}
A real valued process $X$ is \emph{strong Markov} if for every bounded, measurable $g:\reals\rightarrow\reals$ and every $t>0$ there exists a measurable $f:\halfplane\rightarrow\reals$ such that
\begin{equation*}
f(\tau,X_\tau) = \E{g(X_{\tau+t})\mid \setsF_\tau}
\end{equation*}
for every finite stopping time $\tau$.
\end{definition}
We now state the uniqueness result for continuous martingales.
\begin{theorem}\label{thm:cts mart uniqueness}
Let $X$ and $Y$ be continuous and strong Markov martingales such that, for every $t\in\reals_+$, $X_t$ and $Y_t$ have the same distribution. Then, $X$ and $Y$ have the same joint distribution.
\end{theorem}
The proof that this follows from our generalized backward equation is given in Section \ref{sec:cts martingales}. The question of existence of continuous martingales with given marginals will be studied in a future paper.

Let us briefly discuss the idea here in the case where $X$ is a martingale, so the coefficient $b$ in SDE (\ref{eqn:SDE}) is zero. Then, let $C:\halfplane\rightarrow\reals$ be the function
\begin{equation}\label{eqn:C def}
C(t,x)=\E{(X_t-x)_+}.
\end{equation}
Clearly, $C(t,x)$ depends only on the distribution of $X_t$. Furthermore, it is easy to see that knowing $C$ is equivalent to knowing all the marginal distributions of $X$. If $X$ has smooth probability densities, the following follows from the Fokker-Planck (or forward Kolmogorov) equation
\begin{equation}\label{eqn:C pde}
\frac{\partial}{\partial t}C(t,x)=\frac{1}{2}\sigma^2(t,x)\frac{\partial^2}{\partial x^2}C(t,x).
\end{equation}
This equation is well known in finance, where it is used to calibrate the local volatility model to observed option prices (see \citep{Dupire2}).
We note that it is possible to combine this with (\ref{eqn:Feynman Kac}) in order to eliminate $\sigma$, giving the following martingale condition for $f(t,X_t)$,
\begin{equation}\label{eqn:Smooth Martingale PDE:C Martingale condition}
\frac{\partial f}{\partial t}\frac{\partial^2 C}{\partial x^2} + \frac{\partial C}{\partial t}\frac{\partial^2 f}{\partial x^2}=0.
\end{equation}
Note that, although we have removed any reference to SDE (\ref{eqn:SDE}), the equation is still only defined when $f$ and $C$ are twice continuously differentiable. The idea now is to smooth this expression by multiplying by a twice differentiable function $\theta$ with compact support in $(0,\infty)\times\reals$ and integrating. To this end, make the following definition,
\begin{equation}\label{eqn:mufc smooth 1}
\mu_{[f,C]}(\theta)=\dblint \theta\left(\frac{\partial f}{\partial t}\frac{\partial^2 C}{\partial x^2} + \frac{\partial C}{\partial t}\frac{\partial^2 f}{\partial x^2}\right)\,dt\,dx.
\end{equation}
So $\mu_{[f,C]}(\theta)$ is a linear function of $\theta$, and is symmetric in exchanging $f$ and $C$.
Integration by parts gives
 \begin{equation}\label{eqn:mufc smooth 2}
\mu_{[f,C]}(\theta)= \dblint\left( \frac{\partial f}{\partial x}\frac{\partial C}{\partial x}\frac{\partial \theta}{\partial t}-\frac{\partial \theta}{\partial x}\frac{\partial f}{\partial x}\frac{\partial C}{\partial t} - \frac{\partial C}{\partial x}\frac{\partial \theta}{\partial x}\frac{\partial f}{\partial t}\right)\,dt\,dx
\end{equation}
and the martingale condition for $f(t,X_t)$ can be written as $\mu_{[f,C]}=0$.
In order for the right hand side of (\ref{eqn:mufc smooth 2}) to be well defined, $f$ and $C$ only have to be once differentiable.
In fact, we can go further and remove the first order derivatives with respect to $t$, replacing them with Lebesgue-Stieltjes integrals. For every $x\in\reals$ such that $f(t,x)$ is right-continuous with locally finite variation in $t$, the Lebesgue-Stieltjes integral $\int\cdot\,d_tf(t,x)$ is (locally) a finite signed measure satisfying
\begin{equation*}
\int_{t_0}^{t_1}\,d_tf(t,x) = f(t_1,x)-f(t_0,x).
\end{equation*}
In particular, if $f$ is continuously differentiable in $t$ then $d_tf=(\partial f/\partial t)\,dt$.
Furthermore, for the right hand side of (\ref{eqn:mufc smooth 2}) to be well defined, it is only necessary for the partial derivatives with respect to $x$ to exist in an almost-everywhere sense.
In particular, equation (\ref{eqn:C def}) implies that $C(t,x)$ will be convex in $x$, so its partial derivatives exist almost everywhere, and Jensen's inequality implies that it will be increasing in $t$. So, expression (\ref{eqn:mufc smooth 2}) is well defined even if $C$ is not smooth. Similarly if $f$ is defined by (\ref{eqn:f is exp g}), $g$ is convex, and $X$ is a continuous and strong Markov martingale then $f(t,x)$ will be convex in $x$ and decreasing in $t$ (see \citep{Lowther1}). Then, as we shall see, this leads to a well-defined expression for $\mu_{[f,C]}(\theta)$.

Although the above argument gives us reason to believe that $\mu_{[f,C]}=0$ is a valid martingale condition for $f(t,X_t)$, it only proves it in the case where $f$ and $C$ are smooth, and $X$ satisfies the SDE (\ref{eqn:SDE}).
The general case is much more difficult, because results such as It\^o's formula can't be applied directly to non-differentiable functions.
Proving that $\mu_{[f,C]}=0$ is a sufficient condition for $f(t,X_t)$ to be a martingale is the hardest direction, because we do not even know a~priori that it is a semimartingale. This will require the results of \citep{Lowther3} which show that $f(t,X_t)$ is a \emph{Dirichlet process}, so can be decomposed into local martingale and zero continuous quadratic variation terms. A time reversal argument is applied in Sections \ref{sec:reversed condvar} and \ref{sec:pf of 2.5 and 2.6} to show that the zero continuous quadratic variation term has finite variation.

The analysis above can be extended to more general jump diffusions given by a time dependent generator
\begin{equation}\label{eqn:L def}\begin{split}
L_tf(x) ={}& \frac{1}{2}\sigma^2(t,x)\frac{\partial^2 f(x)}{\partial x^2}+b(t,x)\frac{\partial f(x)}{\partial x}\\
&+\int\left(f(y)-f(x)-(y-x)\frac{\partial f(x)}{\partial x}\right)j(t,x,y)\,dy,
\end{split}\end{equation}
where $b,\sigma,j$ are sufficiently well-behaved functions. Again, it is possible to eliminate $\sigma$ from the forward and backward equations. We do not go through the details here as it is the same in principle to the above argument except that there are additional terms due to the drift $b$ and jump component $j$. The following martingale condition for $f(t,X_t)$ is then obtained,
\begin{equation}\label{eqn:jump diff mart cond}\begin{split}
&\mu_{[f,C]}(\theta)+\int\!\!\!\int\int_{-\infty}^x \theta\pd{2}(t,y)f\pd{2}(t,y)\,dy\,p_t(x)b(t,x)\,dx\,dt\\
&+\int\!\!\!\int\!\!\!\int J_t(y,z)\,dy\,dz\,dt=0
\end{split}\end{equation}
where $p_t(x)=C\pd{22}(t,x)$ is the probability density of $X_t$ and $J_t(y,z)$ is
\begin{equation*}
J_t(y,z)=p_t(y)j(t,y,z)\int_{y}^z \left(f(t,x)-f(t,z)+(z-x)f\pd{2}(t,x)\right)\theta\pd{2}(t,x)\,dx.
\end{equation*}
All of the terms in expression (\ref{eqn:jump diff mart cond}) can be defined for general quasimartingales, as shown in Section \ref{sec:quasimartingales}, and will give the required martingale condition.

\section{Continuous Martingales}
\label{sec:cts martingales}

In this section we state the martingale condition for $f(t,X_t)$ in the case where $X$ is a continuous and strong Markov martingale, and show how it implies that such processes are uniquely determined by their marginal distributions.
As discussed above, the martingale condition will be of the form $\mu_{[f,C]}(\theta)=0$ which, for continuously differentiable functions, is defined by expression (\ref{eqn:mufc smooth 2}).
We will extend this definition to the following class of non-differentiable functions for $\theta$, $f$ and $C$. Here, $\int_0^T \,|d_tf(t,x)|$ represents the variation of $f$ with respect to $t$ over the interval $[0,T]$.

\begin{definition}\label{def:nicefunc}
Denote by $\nicefunc$ the set of functions $f:\halfplane\rightarrow\reals$ such that
\begin{itemize}
\item $f(t,x)$ is Lipschitz continuous in $x$ and \cadlag\ in $t$,
\item for every $K_0<K_1\in\reals$ and $T\in\reals_+$ then
\begin{equation*}
\dlimint{K_0}{K_1}{0}{T}\,|d_tf(t,x)|\,dx<\infty,
\end{equation*}
\item the left and right derivatives of $f(t,x)$ with respect to $x$ exist everywhere.
\end{itemize}
We also define $\nicefuncK$ to consist of those functions $f\in\nicefunc$ with compact support in $(0,\infty)\times\reals$.
\end{definition}

Except for the global Lipschitz condition, which is adopted here for convenience, this  class of functions is the same as that used for the decomposition results in \citep{Lowther3}. Note that if $f(t,x)$ is convex and Lipschitz continuous in $x$, and right continuous and monotonic in $t$ then $f\in\nicefunc$.
In particular, $C(t,x)$ defined by (\ref{eqn:C def}) is in $\nicefunc$ whenever $X$ is a right-continuous martingale.
Furthermore, if $f\in\nicefunc$ then $f(t,x)$ will have locally finite variation in $t$ for almost every $x$, and the integral $\dblint\cdot\,d_tf(t,x)\,dx$ is well defined and locally gives a finite signed measure.

Let us now introduce use a bit of simplifying notation. Given any function $f:\halfplane\rightarrow\reals$ where its arguments occur as dummy variables in an expression, then we shall often omit them. In this case we shall always use $s$, $t$ or $u$ for the first argument and $x$ or $y$ for the second. So, for example $\dblint f\,dx\,dt$ is the same as $\dblint f(t,x)\,dx\,dt$.
We shall use $f^-$ to represent the left limit
\begin{equation*}
f^-(t,x)=f(t-,x)\equiv\left\{
\begin{array}{ll}
\lim_{s\uparrow\uparrow t} f(s,x),&\textrm{if $t>0$},\\
f(0,x),&\textrm{if $t=0$}.
\end{array}
 \right.
\end{equation*}
Furthermore, when we write $f^-\pd{2}(t,x)$ then we shall take this to mean the derivative of $f^-(t,x)$ with respect to $x$.

As was noted earlier, in order for expression (\ref{eqn:mufc smooth 2}) for $\mu_{[f,C]}$ to be well defined, the partial derivatives of $f(t,x)$ and $C(t,x)$ with respect to $x$ are required to exist in an almost-everywhere sense. Fortunately this is indeed the case for functions in $\nicefunc$.
\begin{lemma}\label{lemma:nicefunc ae diff}
If $f,g\in\nicefunc$, then $f(t,x)$ and $f^-(t,x)$ are differentiable in $x$ almost everywhere with respect to the measure $\dblint\cdot\,|d_tg(t,x)|\,dx$.
\end{lemma}
\begin{proof}
By Lemma 3.3 of \citep{Lowther3}, $f(t,x)$ is differentiable in $x$ almost everywhere with respect to the measure $\dblint\cdot\,|d_tg(t,x)|\,dx$. We just need to extend this to $f^-$. Note that as $f(t,x)$ is jointly continuous in $x$ and \cadlag\ in $t$, there can only be countably many times $t$ at which $f^-\not= f$.
So, by countable additivity of the measure $\dblint\cdot\,|d_tg(t,x)|\,dx$, we can restrict to fixed times $T>0$. As $f^-(T,x)$ is Lipschitz continuous in $x$, Lebesgue's theorem says that it will be differentiable almost-everywhere with respect to the Lebesgue measure (see \citep{Heinonen} Theorem 3.2). So,
\begin{equation*}\begin{split}
&\dblint1_{\{t=T\}}1_{\{f^-(t,x)\textrm{ is not differentiable in }x\}}\,|d_tg(t,x)|\,dx\\
={}&\int 1_{\{f^-(T,x)\textrm{ is not differentiable in }x\}}\,|g(T,x)-g(T-,x)|\,dx=0
\end{split}\end{equation*}
as required.
\end{proof}
This allows $\mu_{[f,C]}(\theta)$ to be extended to all $f,C\in\nicefunc$ and $\theta\in\nicefuncK$.
\begin{definition}\label{defn:mufc}
For every $f,g\in\nicefunc$ define the linear map $\mu_{[f,g]}:{\nicefuncK}\rightarrow\reals$,
\begin{equation*}
\mu_{[f,g]}\left( \theta \right)=
\dblint f\pd{2}g\pd{2}\,d_t\theta\,dx
-\dblint \theta^-\pd{2} f^-\pd{2}\,d_tg\,dx-\dblint g^-\pd{2}\theta^-\pd{2}\,d_tf\,dx.
\end{equation*}
\end{definition}
In particular, if $f$ and $g$ are twice continuously differentiable this definition coincides with that given by equation (\ref{eqn:mufc smooth 2}), and $\mu_{[f,g]}$ locally defines a signed measure with
\begin{equation}\label{eqn:smooth dmufg}
d\mu_{[f,g]}=
\left(\frac{\partial f}{\partial t}\frac{\partial^2 g}{\partial x^2} + \frac{\partial g}{\partial t}\frac{\partial^2 f}{\partial x^2}\right)\,dt\,dx.
\end{equation}
Then, the martingale condition will be $\mu_{[f,C]}=0$, where $C$ is defined by (\ref{eqn:C def}). In fact, we will go a bit further than this and show that $\mu_{[f,C]}$ determines the drift component of $f(t,X_t)$.

We now define the drift measure of $f(t,X_t)$.
Here, the terminology \emph{local signed measure} on $\halfplane$ is used to mean a real valued function on the bounded Borel subsets of $\halfplane$, which is a finite signed measure when restricted to any bounded Borel set.

\begin{definition}\label{defn:muf}
Let $X$ be a real valued and adapted \cadlag\ process.
Then, we write $\Dom(X)$ for those functions $f\in\nicefunc$ such that the following decomposition exists,
\begin{equation}\label{eqn:decomposition of f in domain of X}
f(t,X_t) = M_t+A_t.
\end{equation}
Here, $M$ is a \cadlag\ local martingale and $A$ is a \cadlag\ previsible process of locally finite variation with $A_0=0$, and such that $\int 1_{\{(s,X_{s-})\in S\}}\,dA_s$ has integrable variation for every bounded Borel subset $S$ of $\halfplane$.
We then define the local signed measure $\mu^X_f$ by
\begin{equation*}
\mu^X_f(\theta) = \E{\int\theta(t,X_{t-})\,dA_t}
\end{equation*}
for every bounded measurable $\theta:\halfplane\rightarrow\reals$ with bounded support.
\end{definition}

Note that if decomposition (\ref{eqn:decomposition of f in domain of X}) exists then it is unique. If $f(t,X_t)=M_t+A_t=M^\prime_t+A^\prime_t$ were two such decompositions then $M-M^\prime=A^\prime-A$ would be a previsible local martingale with locally finite variation, and so is constant. Therefore $M=M^\prime$ and $A=A^\prime$. So, $\mu^X_f$ is well defined.

We would like to express the martingale property of $f(t,X_t)$ as $\mu^X_f=0$. Certainly, it is clear that this condition will be satisfied whenever $f(t,X_t)$ is a martingale, because then $A=0$ in decomposition (\ref{eqn:decomposition of f in domain of X}). However, the converse is not obvious, and the following result will be required.
Recall that a process $X$ is said to be quasi-left-continuous if $X_{\tau-}=X_\tau$ (a.s.) for all previsible stopping times $\tau>0$.

\begin{theorem}\label{thm:abs drift is abs measure}
Let $X$ be a \cadlag\ and quasi-left-continuous strong Markov process, $f\in\Dom(X)$, and $f(t,X_t)=M_t+A_t$ be decomposition (\ref{eqn:decomposition of f in domain of X}). Then, for every non-negative and measurable $\theta:\halfplane\rightarrow\reals$,
\begin{equation*}
|\mu^X_f|(\theta)=\E{\int\theta(t,X_{t-})\,|dA_t|}.
\end{equation*}
\end{theorem}
The proof of this is left until Section \ref{sec:pf of 2.5 and 2.6}.
Here we are using the notation $|\mu|$ for the variation of a local signed measure $\mu$,
\begin{equation*}
|\mu|(\theta)\equiv\sup_{|g|\le 1}\mu(g\theta)
\end{equation*}
where the supremum is taken over all bounded and measurable $g$ with $|g|\le 1$.
In particular, if $\mu^X_f=0$ then Theorem \ref{thm:abs drift is abs measure} shows that $A$ is constant, so $f(t,X_t)$ is a local martingale.

We now state the main result for functions of continuous martingales, namely that $\mu_{[f,C]}$ describes the drift of $f(t,X_t)$.
\begin{theorem}\label{thm:mufc defines drift of acd process}
Let $X$ be a continuous and strong Markov martingale and define $C\in\nicefunc$ by equation (\ref{eqn:C def}). If $f\in\nicefunc$ then the following are equivalent.
\begin{itemize}
\item $f\in\Dom(X)$.
\item There exists a local signed measure $\mu$ such that $\mu_{[f,C]}(\theta)=\mu(\theta^-)$ for all $\theta\in\nicefuncK$.
\end{itemize}
Furthermore, if these conditions are true then $\mu=\mu^X_f$.
\end{theorem}
The proof of this is again left until Section \ref{sec:pf of 2.5 and 2.6}.
In the remainder of this section, we show how Theorem \ref{thm:mufc defines drift of acd process} implies that continuous and strong Markov martingales are uniquely determined by their marginal distributions. The first step is the promised martingale condition $\mu_{[f,C]}=0$, which is our generalization of the backward Kolmogorov equation (\ref{eqn:Feynman Kac}).

\begin{theorem}\label{thm:acd Martingale Condition}
Let $X$ be a continuous and strong Markov martingale and $f\in\nicefunc$. Then $f(t,X_t)$ is a martingale if and only if $\mu_{[f,C]}(\theta)=0$ for all $\theta\in\nicefuncK$.
\end{theorem}
\begin{proof}
First, suppose that $f(t,X_t)$ is a martingale. Then we can take $A=0$ in decomposition (\ref{eqn:decomposition of f in domain of X}), and it follows that $f\in\Dom(X)$ and $\mu^X_f=0$. Theorem \ref{thm:mufc defines drift of acd process} gives $\mu_{[f,C]}(\theta)=\mu^X_f(\theta^-)=0$.

Conversely, suppose that $\mu_{[f,C]}=0$. Taking $\mu=0$ in Theorem \ref{thm:mufc defines drift of acd process} shows that $f\in\Dom(X)$ and $\mu^X_f=0$.
If $f(t,X_t)=M_t+A_t$ is decomposition (\ref{eqn:decomposition of f in domain of X}) then Theorem \ref{thm:abs drift is abs measure} says that $A$ has zero variation, so $A=0$ and $f(t,X_t)=M_t$ is a local martingale. As $X$ is a martingale and $f(t,x)$ is Lipschitz continuous in $x$, it follows that for every $t>0$, $f(\tau,X_\tau)$ is uniformly integrable over all stopping times $\tau\le t$. So, $f(t,X_t)$ is a proper martingale.
\end{proof}

In order to complete the proof that continuous and strong Markov martingales are uniquely determined by their marginal distributions, we require the following result, which follows from Theorems 1.5 and 1.6 of \citep{Lowther1}.
\begin{theorem}\label{thm:f is convex dec}
Let $X$ be a continuous and strong Markov martingale, $g:\reals\rightarrow\reals$ be convex and Lipschitz continuous, and choose any $T\ge 0$. Then there exists an $f:\halfplane\rightarrow\reals$ such that $f(t,x)$ is convex and Lipschitz continuous in $x$, right-continuous and decreasing in $t$ and
\begin{equation}\label{eqn:f is exp g 2}
f(t,X_t)=\E{g(X_T)|\setsF_t}
\end{equation}
for all $t\le T$.
\end{theorem}
In particular, $f\in\nicefunc$. The proof of Theorem \ref{thm:cts mart uniqueness} now follows easily.

\begin{proof}[Proof of Theorem \ref{thm:cts mart uniqueness}]
Note that for any $T>0$, the stopped processes $X^T$ and $Y^T$ are continuous and strong Markov martingales.
Define $C\in\nicefunc$ by $C(t,x)=\E{(X^T_t-x)_+}$ and choose any convex and Lipschitz continuous $g:\reals\rightarrow\reals$.
Then choose $f\in\nicefunc$ such that equality (\ref{eqn:f is exp g 2}) is satisfied for every $t\le T$.
By definition, $f(t,X^T_t)$ is a martingale, so Theorem \ref{thm:acd Martingale Condition} gives $\mu_{[f,C]}=0$.

As $X$ and $Y$ have the same marginal distributions, $C$ is also given by $C(t,x)=\E{(Y^T_t-x)_+}$. So Theorem \ref{thm:acd Martingale Condition} says that $f(t,Y^T_t)$ is a martingale.
In particular, equality (\ref{eqn:f is exp g 2}) is satisfied with $Y$ in place of $X$, showing that $(X_t,X_T)$ and $(Y_t,Y_T)$ have the same distribution whenever $t<T$.
As $X$ and $Y$ are required to be Markov and have the same initial distributions, they must have the same joint distributions.
\end{proof}

\section{Quasimartingales}
\label{sec:quasimartingales}

In the previous section the martingale condition $\mu_{[f,C]}=0$ was described for functions of continuous and strong Markov martingales, and it was shown that this implies that such processes are uniquely determined by their marginal distributions. We now turn to the more general case where the process $X$ need not be either continuous or a martingale. The idea here is to use a generalization of (\ref{eqn:jump diff mart cond}). Note that this equation involves the drift $b$ and jump rate $j$ of the process, so it will be necessary to restrict to processes for which a generalization of these terms can be defined. For that reason, we look at quasimartingales.

However, the results stated in this section will not be as strong as for the continuous martingale case.
In particular, we will not extend Theorem \ref{thm:cts mart uniqueness} to quasimartingales. This is a subject which will be taken up in more detail in a later paper.

The simple previsible processes are defined to be linear combinations of processes of the form
$A1_{\srbrack T,\infty\srbrack}$
where $T\ge 0$ and $A$ is an $\setsF_T$-measurable random variable. Then, the definition of quasimartingales used in this paper is as follows.
\begin{definition}\label{def:VarX}\label{def:quasimartingale}
Let $X$ be a real valued stochastic process such that $\E{|X_t|}<\infty$ for every $t\in\reals_+$. Then, define
\begin{equation*}
\VarX{t}\equiv\sup\left\{\E{\int_0^t\xi\,dX}:\xi\textrm{ is simple previsible and }|\xi|\le 1\right\}.
\end{equation*}
Then, a quasimartingale $X$ is a \cadlag\ adapted process such that $\E{|X_t|}<\infty$ and $\VarX{t}<\infty$ for every $t\in\reals_+$.
\end{definition}
This definition of quasimartingales is slightly more general than that used in \citep{HeWangYan}, \citep{Protter} and \citep{Rogers}. This is because we are only really interested in the properties of processes over finite time horizons, and Definition \ref{def:quasimartingale} is the same as saying that the stopped process $X^t$ is a quasimartingale according to the usual definition for every finite time $t$.

We also require a uniform integrability property for the processes. If $\{X_\tau:\tau\textrm{ is a stopping time}\}$ is uniformly integrable, then the process $X$ is said to be of class (D). As we are only interested in finite time horizons, the following definition will be used.

\begin{definition}
A process $X$ is of class (DL) if the stopped process $X^t$ is of class (D) for every $t\in\reals_+$.
\end{definition}

If $C(t,x)$ is defined by equation (\ref{eqn:C def}), then it is in $\nicefunc$ whenever $X$ is a class (DL) quasimartingale. So, $\mu_{[f,C]}$ defined in the previous section is also defined for quasimartingales.

\begin{lemma}\label{lemma:C in nicefunc}
Let $X$ be a quasimartingale of class (DL) and set $C(t,x)=\E{(X_t-x)_+}$. Then, $C(t,x)+\VarX{t}$ is an increasing function of $t$ and $C\in\nicefunc$.
\end{lemma}
\begin{proof}
If $t>s$ then Jensen's inequality gives
\begin{equation*}\begin{split}
C(t,x) &= \E{(X_t-x)_+}\ge \E{(\E{X_t|\setsF_s}-x)_+}\\
&\ge \E{(X_s-x)_+}-\E{\left|\E{X_t-X_s|\setsF_s}\right|}\\
&\ge C(s,x)-\VarX{t}+\VarX{s}
\end{split}\end{equation*}
Therefore, $C(t,x)+\VarX{t}$ is increasing in $t$. So,
\begin{equation*}
\int_{K_0}^{K_1}\int_0^T\,|d_tC(t,x)|\,dx\le \int_{K_0}^{K_1}\left( C(T,x)+2\VarX{T}\right)\,dx<\infty.
\end{equation*}
As $(X_t-x)_+$ is convex in $x$ it follows that $C(t,x)$ is also convex in $x$.
Finally, if we choose $t_n\downarrow t$ then $(X_{t_n}-x)_+$ is a uniformly integrable sequence converging to $(X_t-x)_+$. So, $C(t_n,x)\rightarrow C(t,x)$. This shows that $C(t,x)$ is right-continuous in $t$, so $C\in\nicefunc$.
\end{proof}
It can be shown that this result remains true even without the class (DL) assumption, although that is not required here. Restricting to quasimartingales allows us  split $X$ up into a martingale part and a drift component, using the following version of Rao's decomposition.

\begin{lemma}\label{lemma:Rao decomposition}
If $X$ is a quasimartingale then it has the unique decomposition
\begin{equation}\label{eqn:thm:quasi decomp}
X=M+A
\end{equation}
where $M$ is a local martingale and $A$ is a \cadlag\ previsible process with locally finite variation and $A_0=0$. Furthermore, for every $t>0$,
\begin{equation*}
\E{\int_0^t\,|dA|}\le\VarX{t}
\end{equation*}
and, if $X$ is of class (DL), $M$ is a martingale.
\end{lemma}
\begin{proof}
First, Rao's Theorem (\citep{Protter} III Theorem 15) says that decomposition (\ref{eqn:thm:quasi decomp}) exists and is unique.

Choose any $T>0$ and  consider the quasimartingale $Y_t=X^T_t-\E{X_T|\setsF_t}$. As $Y_t=0$ for $t\ge T$, the definition of ${\rm Var}_X(T)$ above agrees with the conditional variation ${\rm Var}(Y)$ used in \citep{HeWangYan} VIII Section 2. So Rao's decomposition (\citep{HeWangYan} Theorem 8.13) gives $Y=Z^1-Z^2$ for non-negative supermartingales $Z^1,Z^2$ satisfying $\E{Z^1_0+Z^2_0}=\VarX{T}$.
By the Doob-Meyer decomposition (\citep{Protter} III Theorem 13), $Z^i=N^i-B^i$ ($i=1,2$) where $N^i$ are local martingales and $B^i$ are increasing previsible processes with $B^i_0=0$.
So, the uniqueness of decomposition (\ref{eqn:thm:quasi decomp}) gives $A^T=B^2-B^1$ and
\begin{equation*}
\E{\int_0^T\,|dA|}\le\E{B^1_T+B^2_T}\le\E{Z^1_0+Z^2_0}=\VarX{T}.
\end{equation*}
Finally, if $X$ is of class (DL) then so is $M=X-A$, in which case $M$ is a martingale.
\end{proof}

Decomposition (\ref{eqn:thm:quasi decomp}) allows the drift measure of $X$ to be defined as follows.

\begin{definition}
Let $X$ be a quasimartingale and $X=M+A$ be decomposition (\ref{eqn:thm:quasi decomp}). Then, define the \emph{drift measure} of $X$ to be the local signed measure $\mu_X$ given by
\begin{equation*}
\mu_X(\theta) = \E{\int_0^\infty\theta(t,X_{t-})\,dA_t}
\end{equation*}
for all bounded and measurable $\theta:\halfplane\rightarrow\reals$ with bounded support.
\end{definition}

The drift measure $\mu_X$ enables equation (\ref{eqn:jump diff mart cond}) to be extended to more general processes simply by replacing the $p_t(x)b(t,x)\,dx\,dt$ term by $d\mu_X(t,x)$. It remains to extend the jump component $J_t(y,z)$, which is done as follows.

\begin{definition}\label{def:J}
Let $X$ be a \cadlag\ process. For $\theta,f\in\nicefunc$ define
\begin{equation}\label{eqn:def:J}
J^X_t(\theta,f)\equiv\int_{X_{t-}}^{X_t}\left(f(t,x)-f(t,X_t)+(X_t-x)f^-\pd{2}(t,x)\right)\theta^-\pd{2}(t,x)\,dx.
\end{equation}
\end{definition}
This allows equation (\ref{eqn:jump diff mart cond}) to be rewritten in a form that makes sense for all class (DL) quasimartingales.
In order for this expression to be well defined, it is necessary for the sum of the jump terms $J^X_t(\theta,f)$ over all times $t$ to be integrable. This is the case, as will be shown in Lemma \ref{lemma:sum J is integrable}. So, generalizing the left hand side of (\ref{eqn:jump diff mart cond}) gives the following definition.

\begin{definition}\label{def:mutfx}
Let $X$ be a class (DL) quasimartingale. Then, for every $f\in\nicefunc$ define the linear map $\tilde\mu^X_f :{\nicefuncK}\rightarrow\reals$,
\begin{equation*}
\tilde\mu^X_f\left( \theta \right)= \mu_{[f,C]}(\theta) + \int\int_{-\infty}^x\theta\pd{2}^-(t,y)f\pd{2}^-(t,y)\,dy\,d\mu_X(t,x) + \E{\sum_{t>0}J^X_t(\theta,f)}
\end{equation*}
where $C\in\nicefunc$ is defined by $C(t,x)=\E{(X_t-x)_+}$.
\end{definition}

In particular, if $X$ is a continuous martingale then the drift and jump terms above are zero, so $\tilde\mu^X_f=\mu_{[f,C]}$.
Using this definition, we now state the main result of this paper for quasimartingales. Recall that $\mu^X_f$ is the drift measure of $f(t,X_t)$ introduced in Definition \ref{defn:muf}. Note that here we make no requirement for $X$ to be a Markov process.

\begin{theorem}\label{thm:mux for gen proc}
Let $X$ be a class (DL) quasimartingale and $f\in\Dom(X)$. Then $\mu^X_f(\theta^-)=\tilde\mu^X_f(\theta)$ for all $\theta\in\nicefuncK$.
\end{theorem}

The proof of this is left until Section \ref{sec:pf of 4.4}.
In particular, if $f(t,X_t)$ is a martingale then it follows that $f$ is in $\Dom(X)$ and $\mu^X_f=0$. So Theorem \ref{thm:mux for gen proc} gives $\tilde\mu^X_f=0$, and represents a general form of the backward equation (\ref{eqn:Feynman Kac}). Conversely, if $f\in\Dom(X)$ and $\tilde\mu^X_f=0$ this shows that $\mu^X_f=0$ and, if $X$ is also strong Markov and quasi-left-continuous, Theorem \ref{thm:abs drift is abs measure} implies that $f(t,X_t)$ is a martingale.
However, if we only know that $\tilde\mu^X_f=0$ then this does not allow us to infer that $f$ is in $\Dom(X)$, so Theorem \ref{thm:mux for gen proc} cannot be applied.

The reason that we do not obtain as strong results for quasimartingales as in the continuous martingale case is that, for conditional expectations of convex functions of $X$, Theorem \ref{thm:f is convex dec} allowed us to restrict attention to $f\in\nicefunc$. This is required in the proofs of Theorems \ref{thm:cts mart uniqueness} and \ref{thm:mufc defines drift of acd process}. Unfortunately, this result applies only to continuous martingales. In a later paper we will look at conditions on the process $X$ under which $f$ can be approximated sufficiently well by functions in $\Dom(X)$, in order to extend the results for continuous martingales to processes with non-zero drift.

\section{Calculating \texorpdfstring{$\mu^X_f$}{µ} For Smooth Functions}
\label{sec:mufx for smooth f}

We start the proof of Theorem \ref{thm:mux for gen proc} by looking at the case where $f$ is twice continuously differentiable. This is much simpler than for general functions in $\nicefunc$, because It\^o's formula can be applied to explicitly decompose $f(t,X_t)$ into local martingale and finite variation terms.

First, as the quadratic covariation $[X,Y]$ of semimartingales $X$ and $Y$ is a \cadlag\ process with locally finite variation whose jumps satisfy $\Delta[X,Y]=\Delta X\Delta Y$, its continuous part can be written as
\begin{equation}\label{eqn:cts quad cov}
[X,Y]^c_t=[X,Y]_t - \sum_{s\le t}\Delta X_s\Delta Y_s.
\end{equation}
We also write $[X]^c\equiv[X,X]^c$ for the the continuous part of the quadratic variation of $X$.
Then, the following lemma will be used to evaluate terms of the form $\E{\int\theta(t,X_t)\,d[X]^c}$.

\begin{lemma}\label{lemma:vol measure}
Let $X$ be a quasimartingale of class (DL). Define $C\in\nicefunc$ by $C(t,x)=\E{(X_t-x)_+}$ and let $\mu_X$ be the drift measure of $X$.
Then, for any bounded and measurable $\theta:\halfplane\rightarrow\reals$ with bounded support,
\begin{align}
&\E{\int_0^\infty\left|\theta(t,X_{t})\right|\,d[X]^c_t}<\infty,\label{eqn:lemma:vol measure:2}\\
&\E{\sum_{t>0}\int_{X_{t-}}^{X_t}(X_t-x)\left|\theta(t,x)\right|\,dx}<\infty,\label{eqn:lemma:vol measure:3}
\end{align}
and
\begin{equation}\label{eqn:lemma:vol measure:1}\begin{split}
\dblint\theta\,d_tC\,dx={}&\frac{1}{2}\E{\int_0^\infty\theta(t,X_{t})\,d[X]^c_t}
+\int\int_{-\infty}^x\theta(t,y)\,dy\,d\mu_X(t,x)\\
&+ \E{\sum_{t>0}\int_{X_{t-}}^{X_t}(X_t-x)\theta(t,x)\,dx}.
\end{split}\end{equation}
\end{lemma}
\begin{proof}
By linearity, it is enough to prove this for non-negative $\theta$, for which equation (\ref{eqn:lemma:vol measure:1}) follows from substituting the definition of $\mu_X$ into Lemma 4.3 of \citep{Lowther3}. As $C\in\nicefunc$ (Lemma \ref{lemma:C in nicefunc}) and $\mu_X$ is finite on finite time intervals, the left hand side and the second term on the right hand side of (\ref{eqn:lemma:vol measure:1}) are finite. Then, as the remaining terms on the right hand side are non-negative, they must also be finite.
\end{proof}

Inequality (\ref{eqn:lemma:vol measure:3}) shows that the sum of the jump terms $J_t(\theta,f)$ used in Definition \ref{def:mutfx} is integrable, so $\tilde\mu^X_f$ is well defined.

\begin{lemma}\label{lemma:sum J is integrable}
Let $X$ be a class (DL) quasimartingale, $f\in\nicefunc$ and $\theta\in\nicefuncK$. Then
$\E{\sum_{t>0}|J^X_t(\theta,f)|}<\infty.
$\end{lemma}
\begin{proof}
If $|f\pd{2}|\le L$,
\begin{equation*}\begin{split}
J_t^X(\theta,f) &= \int_{X_{t-}}^{X_t}\left((f(t,x)-f(t,X_t)+(X_t-x)f\pd{2}^-(t,x)\right)\theta^-\pd{2}(t,x)\,dx,\\
&\le \int_{X_{t-}}^{X_t}(X_t-x)(L+|f\pd{2}^-|)\theta^-\pd{2}\,dx,
\end{split}\end{equation*}
and inequality (\ref{eqn:lemma:vol measure:3}) gives the result.
\end{proof}

Lemma \ref{lemma:vol measure} can be used together with It\^o's formula to prove Theorem \ref{thm:mux for gen proc} when $f$ is continuously differentiable, as we now show. Here we restrict to functions $f$ with compact support. This is just for simplicity, and in any case it will be extended to general $f\in\Dom(X)$ later.

\begin{lemma}\label{lemma:muf is mutf for smooth f}
Let $X$ be a class (DL) quasimartingale and $f\in\nicefuncK$ be twice continuously differentiable. Then $f\in\Dom(X)$ and
$\mu^X_f(\theta^-)= \tilde\mu^X_f(\theta)$
for all $\theta\in\nicefuncK$.
\end{lemma}
\begin{proof}
Let $X=M+A$ be decomposition (\ref{eqn:thm:quasi decomp}) and apply It\^o's formula,
\begin{align*}
& f(T,X_T) = N_T+V_T,\\
& N_T = f(0,X_0)+\int_0^Tf\pd{2}\,dM_t,\\
& V_T = \frac{1}{2}\int_0^T f\pd{22}\,d[X]_t^c+\int_0^T f\pd{2}\,dA_t+\int_0^T f\pd{1}\,dt+\sum_{t\le T}J^0_t,\\
& J^0_t=f(t,X_t)-f - f\pd{2}\Delta X_t=\int_{X_{t-}}^{X_t}(X_t-x)f\pd{22}(t,x)\,dx.
\end{align*}
Here, for brevity, we have omitted the arguments for all functions whenever they are $(t,X_{t-})$.
As $A$ has integrable variation over finite time intervals, inequalities (\ref{eqn:lemma:vol measure:2}) and (\ref{eqn:lemma:vol measure:3}) show that $V$ has integrable variation, so $f$ is in $\Dom(X)$.
It follows from the expression above for $V$ that, for any $\theta\in\nicefuncK$,
\begin{equation*}
\mu^X_f(\theta^-) = \frac{1}{2}\E{\int\theta^- f\pd{22}\,d[X]_t^c}+\int\E{\theta^- f\pd{1}}dt+\mu_X(\theta^- f\pd{2})+\E{\sum_{t>0}\theta^- J^0_t}.
\end{equation*}
By substituting $\theta^-f\pd{22}$ in place of $\theta$ in equation (\ref{eqn:lemma:vol measure:1}), we can eliminate $[X]^c$ from the above to get
\begin{equation*}
\mu^X_f(\theta^-)=P+Q+\E{\sum_{t>0}J^1_t},
\end{equation*}
where $P,Q$ and $J^1_t$ are as follows,
\begin{align*}
P &= \dblint \theta^-f\pd{22}\,d_tC\,dx+\int\E{\theta^- f\pd{1}}\,dt\\
Q &=\mu_X(\theta^- f\pd{2})-\int\int_{-\infty}^x\theta^-(t,y)f\pd{22}(t,y)\,dy\,d\mu_X(t,x)\\
&= \int\int_{-\infty}^x\theta\pd{2}^-(t,y)f\pd{2}(t,y)\,dy\,d\mu_X(t,x)\\
J^1_t &= \theta^-(t,X_{t-})J^0_t -\int_{X_{t-}}^{X_t}(X_t-x)\theta^-(t,x)f\pd{22}(t,x)\,dx\\
&=\int_{X_{t-}}^{X_t}(X_t-x)(\theta^-(t,X_{t-})-\theta^-(t,x))f\pd{22}(t,x)\,dx
\end{align*}
Here, we have applied an integration by parts to $Q$, and applying a further couple of integration by parts to $J^1_t$ gives $J^1_t=J^X_t(\theta,f)$.
It only remains to be shown that $P=\mu_{[f,C]}(\theta)$. We use the following identity which applies to all $g,h\in\nicefuncK$
\begin{equation*}
\dblint g^-\pd{2}\,d_th\,dx = \dblint h\pd{2}\,d_tg\,dx.
\end{equation*}
If either $g$ or $h$ is twice continuously differentiable then this follows from integration by parts. More generally, it is given by \citep{Lowther3}, Lemma 3.1.
Substitute $\theta f\pd{2}$ in place of $g$ and $C$ in place of $h$, and again apply integration by parts,
\begin{equation*}\begin{split}
& \dblint (\theta^-f\pd{22}+\theta\pd{2}^-f\pd{2})\,d_tC\,dx
=\dblint C\pd{2} d_t\left(\theta f\pd{2}\right)\,dx\\
&= \dblint f\pd{2}C\pd{2}\,d_t\theta\,dx+\dblint \theta C\pd{2}f\pd{21}\,dt\,dx\\
&= \dblint f\pd{2}C\pd{2}\,d_t\theta\,dx-\dblint \theta\pd{2} C\pd{2}f\pd{1}\,dt\,dx
-\dblint\theta f\pd{1}\,d_xC\pd{2}\,dt\,dx.
\end{split}\end{equation*}
So, as required, Definition \ref{defn:mufc} of $\mu_{[f,C]}(\theta)$ gives
\begin{equation*}
\mu_{[f,C]}(\theta)=\dblint\theta^-f\pd{22}\,d_tC\,dx+\dblint\theta^-f\pd{1}\,d_xC\pd{2}\,dt=P.\qedhere
\end{equation*}
\end{proof}

\section{Dirichlet Processes}
\label{sec:gen drift}

As we saw in the previous section, the identity $\mu^X_f(\theta^-)=\tilde\mu^X_f(\theta)$ follows from It\^o's formula when $f$ is twice continuously differentiable. However, this method cannot be applied to non-differentiable functions. Instead, we use the results of \citep{Lowther3}, where it was shown that $f(t,X_t)$ can be decomposed into a local martingale and a term whose quadratic variation exists and has zero continuous part.
We now look at some of the properties of these \emph{Dirichlet processes}, which will be required in the remainder of this paper. Such processes form a generalization of semimartingales, and allow decomposition (\ref{eqn:decomposition of f in domain of X}) to be generalized to all $f\in\nicefunc$.

We start with the definition of the quadratic variation for general \cadlag\ processes.
A (stochastic) partition $P$ of $\reals_+$ is a sequence of stopping times $0=\tau_0\le \tau_1\le \tau_2\le\cdots$ such that $\tau_n\rightarrow\infty$ as $n\rightarrow\infty$. For any partition $P$ we write the sequence it defines as $\tau_n^P$, so
\begin{equation*}
0=\tau^P_0\le \tau^P_1\le \tau^P_2\le\cdots\uparrow\infty.
\end{equation*}
Use $|P|$ to denote the maximum step size of $P$,
\begin{equation*}
|P|=\sup_{k\in\nat}\Vert\tau^P_k-\tau^P_{k-1}\Vert_\infty.
\end{equation*}
Here $\Vert U\Vert_\infty$ is the smallest positive real number that is almost surely an upper bound for $|U|$.
Then, for \cadlag\ processes $X,Y$ the approximation $[X,Y]^P$ to the quadratic covariation along a partition $P$ is
\begin{equation}\label{eqn:qv approx def}
[X,Y]^P_t\equiv\sum_{k=1}^\infty(X_{t^P_{k}\wedge t}-X_{t^P_{k-1}\wedge t})(Y_{t^P_{k}\wedge t}-Y_{t^P_{k-1}\wedge t}).
\end{equation}
The quadratic covariation is then defined as the limit of $[X,Y]^P$ as the partition step size $|P|$ goes to zero. This limit is taken in the topology of uniform convergence on compacts in probability (ucp for short).  
\begin{definition}\label{defn:covariation}
Let $X$, $Y$ be real valued \cadlag\ processes. Then their quadratic covariation $[X,Y]$, if it exists, is a \cadlag\ process such that
\begin{equation*}
[X,Y]^P\ucp [X,Y]
\end{equation*}
as $|P|\rightarrow 0$.
Given any \cadlag\ process $X$ then its quadratic variation $[X]$, if it exists, is defined to be $[X]\equiv[X,X]$.
\end{definition}

From uniform convergence, the jumps of $[X,Y]$ must satisfy $\Delta[X,Y]=\Delta X\Delta Y$.
Whenever $[X],[Y],[X,Y]$ all exist then $[X,Y]=[X+Y]/2-([X]+[Y])/2$ is a difference of increasing processes and so has locally finite variation. As with the semimartingale case, the continuous part of the covariation, $[X,Y]^c$, is defined by (\ref{eqn:cts quad cov}), and $[X]^c\equiv[X,X]^c$.

We shall say that $V$ is a \emph{zero continuous quadratic variation} process (\zcqv\ process for short) if it is \cadlag, adapted and its quadratic variation exists and satisfies $[V]^c=0$.
In particular, every \cadlag\ locally finite variation process has zero continuous quadratic variation (see, for example, \citep{Protter} II Theorem 26).
The following result was proven in \citep{Lowther3} Lemma 1.5.
\begin{lemma}\label{lemma:cov with zero qv proc}
Let $X$ be a \cadlag\ process whose quadratic variation exists and $V$ be a \zcqv\ process.
Then, the covariation $[X,V]$ exists and satisfies
\begin{equation*}
[X,V]_t=\sum_{s\le t}\Delta X_s\Delta V_s.
\end{equation*}
\end{lemma}
This is equivalent to saying that $[X,V]^c=0$ and it follows that the class of \zcqv\ processes is closed under taking linear combinations. Furthermore, if $A$ is a previsible \zcqv\ process and $X$ is a semimartingale,
\begin{equation}\label{eq:cov of X and A}
[A,X]_t=\sum_{s\le t}\Delta A_s\Delta X_s=\int_0^t\Delta A\,dX.
\end{equation}
In the special case where $X$ is a local martingale, this shows that $[A,X]$ will also be a local martingale. We use the following definition.

\begin{definition}
A process $X$ is a \emph{Dirichlet process} iff it has the decomposition $X=M+V$ for a local martingale $M$ and a \zcqv\ process $V$.
\end{definition}

Such processes were introduced by Follmer in \citep{Follmer2}, although only continuous processes were considered. The definition was extended to non-continuous processes in \citep{Stricker1} and \citep{Coquet}, where a Dirichlet process is defined as the sum of a semimartingale and a process with zero quadratic variation. Here we are using the more general definition from \citep{Lowther3} where the zero continuous quadratic variation component is also allowed to be non-continuous.
In particular, a semimartingale is a sum of a local martingale and a locally finite variation process, so is a Dirichlet process. Note that Lemma \ref{lemma:cov with zero qv proc} implies that the quadratic covariation between any two Dirichlet processes is well defined.
The reason for our interest in Dirichlet processes is the following decomposition.
 \begin{lemma}\label{lemma:f is M plus A for nondif}
 Let $X$ be a semimartingale such that $X^*_t\equiv \sup_{s\le t}|X_s|$ is locally integrable, and $f\in\nicefunc$. Then, there exists a unique decomposition
 \begin{equation}\label{eqn:lemma:f is M plus A for nondif}
 f(t,X_t)=M_t+A_t
 \end{equation}
where $M$ is a local martingale and $A$ is a previsible \zcqv\ process with $A_0=0$.
 \end{lemma}
 \begin{proof}
 First, by Theorem 1.4 of \citep{Lowther3}, $f(t,X_t)$ is a Dirichlet process. As $f(t,x)$ is Lipschitz continuous in $x$, $\sup_{s\le t}|f(s,X_s)|$ is locally integrable, so decomposition (\ref{eqn:lemma:f is M plus A for nondif}) follows from Lemma 1.6 of \citep{Lowther3}.
\end{proof}
As every \cadlag\ process with locally finite variation also has zero continuous quadratic variation, decompositions (\ref{eqn:decomposition of f in domain of X}) and (\ref{eqn:lemma:f is M plus A for nondif}) coincide in the case that $f\in\Dom(X)$.
We shall use the process $A$ in decomposition (\ref{eqn:lemma:f is M plus A for nondif}) as a generalization of the drift of $f(t,X_t)$.

It will be necessary to integrate with respect to to Dirichlet processes. It is possible to do this by taking limits over partitions in a similar way to the definition of the quadratic covariation above. However, whatever definition is used, by analogy with the case where $Y$ is a semimartingale the following integration by parts formula would be expected to hold.
For any semimartingale $X$ and adapted Dirichlet process $Y$,
\begin{equation}\label{eq:int by parts for gen int}
X_tY_t=X_0Y_0+\int_0^tX_{s-}\,dY_s+\int_0^tY_{s-}\,dX_s+[X,Y]_t.
\end{equation}
This can be taken as the definition of $\int_0^tX_{s-}\,dY_s$, as all the remaining terms are already defined.
Then, the jumps of $U\equiv\int X_{s-}\,dY_s$ satisfy
\begin{equation*}
\Delta U_t=X_{t-}\Delta Y_t,
\end{equation*}
which follows directly from (\ref{eq:int by parts for gen int}).
In the case where $A$ is a previsible \zcqv\ process and $X$ is a semimartingale, combining equations (\ref{eq:int by parts for gen int}) and (\ref{eq:cov of X and A}) gives
\begin{equation}\label{eqn:int wrt prev zcqv}
X_tA_t=X_0A_0+\int_0^t X_{s-}\,dA_s+\int_0^t A_s\,dX_s.
\end{equation}
Furthermore, the jumps of $B=\int X_{s-}\,dA_s$ are $\Delta B_t=X_{t-}\Delta A_t$, which is previsible. So $B$ will also be previsible. In fact, it will also have zero continuous quadratic variation as the following result shows.
\begin{lemma}\label{lemma:int of zcqv is zcqv}
Let $X$ be a semimartingale and $V$ be a \zcqv\ process. Then, $\int X_{s-}\,dV_s$ also has zero continuous quadratic variation.
\end{lemma}
\begin{proof}
We make use of the result that if $f:\reals\rightarrow\reals$ is continuously differentiable and $Z=X+V$ then
\begin{equation*}
f(Z_t)-\int_0^tf^\prime(Z_{s-})\,dX_s
\end{equation*}
has zero continuous quadratic variation. This is a special case of Theorem 2.1 of \citep{Lowther3}. Applying this result with $f(x)=x^2/2$ shows that
\begin{equation}\label{eq:pf:int of zcqv is zcqv:1}
\frac{1}{2}(X_t+V_t)^2-\int_0^t(X_{s-}+V_{s-})\,dX_s
\end{equation}
has zero continuous quadratic variation. Substituting $X=0$ and then $V=0$ into (\ref{eq:pf:int of zcqv is zcqv:1}) shows that $V_t^2/2$ and $X_t^2/2-\int_0^t X_{s-}\,dX_s$ are also \zcqv\ processes. Subtract these from (\ref{eq:pf:int of zcqv is zcqv:1}) to get the \zcqv\ process
\begin{equation*}
X_tV_t-\int_0^t V_{s-}\,dX_s = \int_0^t X_{s-}\,dV_s + X_0V_0+[X,V]_t.
\end{equation*}
As $X_0V_0+[X,V]_t$ locally has finite variation, it follows that $\int X_{s-}\,dV_s$ has zero continuous quadratic variation.
\end{proof}
The associativity of stochastic integration also extends to integration with respect to Dirichlet processes.
\begin{lemma}
Let $X$ be a Dirichlet process, $Y,Z$ be semimartingales, and $U_t\equiv\int_0^t Y_{s-}\,dX_s$. Then,
\begin{equation}\label{eq:lemma:int of zcqv is zcqv:1}
\int_0^t Z_{s-}\,dU_s = \int_0^t Z_{s-}Y_{s-}\,dX_s.
\end{equation}
\end{lemma}
\begin{proof}
If $X$ is a semimartingale, identity (\ref{eq:lemma:int of zcqv is zcqv:1}) follows from standard rules of stochastic calculus. So, by linearity, we may restrict to the case where $X$ is a \zcqv\ process with $X_0=0$. For brevity, define $X\dotp Y$ to be the integral
\begin{equation*}
(X\dotp Y)_t \equiv \int_0^t X_{s-}\,dY_s,
\end{equation*}
which is well defined whenever $X,Y$ are semimartingales, or when one is a semimartingale and the other is an adapted \zcqv\ process. Then, by the normal rules of stochastic integration, $X\dotp (Y\dotp Z)=(XY)\dotp Z$ whenever $Z$ is a semimartingale.
We will also write $X\sim Y$ to mean that $X-Y$ is a locally finite variation process with zero continuous part (i.e., a pure jump process). So $X=Y$ if and only if $X\sim Y$, $\Delta X=\Delta Y$ and $X_0=Y_0$. Then,
\begin{equation*}\begin{split}
Z\dotp (Y\dotp X) &\sim Z(Y\dotp X) - (Y\dotp X)\dotp Z\\
&= ZYX - Z(X\dotp Y+[X,Y])-(Y\dotp X)\dotp Z\\
&\sim ZYX-(ZX)\dotp Y - (X\dotp Y+[X,Y])\dotp Z-[Z,X\dotp Y]-(Y\dotp X)\dotp Z\\
&=ZYX - (XZ)\dotp Y-(XY)\dotp Z-X\dotp [Z,Y]\\
&= ZYX - X\dotp (ZY)\\
&\sim (ZY)\dotp X.
\end{split}\end{equation*}
Here we have used one integration by parts per line, and dropped any term that is just a quadratic covariation involving $X$ or $Y\dotp X$, as these will be pure jump processes. The terms involving $X_0$ have also been dropped because we are restricting to $X_0=0$.
Finally,
\begin{equation*}
\Delta(Z\dotp(Y\dotp X))=Z_-\Delta(Y\dotp X)=Z_-Y_-\Delta X=\Delta((ZY)\dotp X),
\end{equation*}
so $Z\dotp (Y\dotp X)=(ZY)\dotp X$ as required.
\end{proof}

This allows Definition \ref{defn:muf} of $\mu^X_f$ to be extended to general $f\in\nicefunc$.
\begin{definition}\label{def:generalized muf}
Let $X$ be a class (DL) quasimartingale, $f\in\nicefunc$, and let $f(t,X_t)=M_t+A_t$ be decomposition (\ref{eqn:lemma:f is M plus A for nondif}). If $\theta\in\nicefuncK$ is such that $\theta(t,X_t)$ is a semimartingale, then we say that $\theta^-$ is $\mu^X_f$ integrable if and only if
\begin{equation*}
\E{\sup_{t\ge 0}\left|\int_0^t \theta^-(s,X_{s-})\,dA_s\right|}<\infty.
\end{equation*}
In that case we define
$
\mu^X_f( \theta^-)\equiv\E{\int \theta^-(s,X_{s-})\,dA_s}.
$
\end{definition}

\section{Proof of Theorem \ref{thm:mux for gen proc}}
\label{sec:pf of 4.4}

In Lemma \ref{lemma:muf is mutf for smooth f} it was shown that $\mu^X_f(\theta^-)=\tilde\mu^X_f(\theta)$ for twice continuously differentiable functions $f$. The method that we shall use to extend this result to more general functions is to employ integration by parts to exchange the roles of $\theta$ and $f$. The idea is to show that $\mu^X_f$ and $\tilde\mu^X_f$ both satisfy the same integration by parts formula, which will prove the result when $f$ has compact support. The extension to general $f\in\nicefunc$ will make use of the following result. Note that here Definition \ref{def:generalized muf} for $\mu^X_f$ is being employed, which does not require $f$ to be in $\Dom(X)$.

\begin{lemma}\label{lemma:mufx for disjoint supports}
Let $X$ be a class (DL) quasimartingale, $f\in\nicefunc$ and $\theta\in\nicefuncK$ be such that $\theta(t,X_t)$ is a semimartingale. If $f$ and $\theta$ have disjoint supports then $\theta^-$ is $\mu^X_f$-integrable and $\mu^X_f(\theta^-)=\tilde\mu^X_f(\theta)$.
\end{lemma}
\begin{proof}
Define the increasing sequences of stopping times $\sigma_n,\tau_n$ by $\tau_0=0$ and
\begin{equation*}
\sigma_{n} = \inf\left\{ t\ge \tau_{n-1}:\theta(s,X_s)\not= 0\right\},\ \tau_{n} = \inf\left\{ t\ge \sigma_n:f(s,X_s)\not=0\right\}
\end{equation*}
for $n\in\nat$.
Note that $(\sigma_n,X_{\sigma_n})$ is in the support of $f$ and $(\tau_n,X_{\tau_n})$ is in the support of $\theta$. If these sequences converge to a finite limit $\tau$, then $(\tau,X_{\tau-})$ would be in the intersection of the supports of $f$ and $\theta$. However, this intersection is empty, so $\tau_n$ and $\sigma_n$ increase to infinity.

Set $Y_t=f(t,X_t)$. Then, the pure jump process $V\equiv\sum_{n=1}^\infty 1_{\slbrack\sigma_n,\tau_n\slbrack}$ satisfies $VY=0$, so integration by parts gives
\begin{equation*}
\int_0^tV_{s-}\,dY_s
=
[V,Y]_t-\int_0^t Y_{s-}\,dV_s
=\sum_{s\le t}(\Delta V_s\Delta Y_s-Y_{s-}\Delta V_s)=\sum_{s\le t}V_{s-}Y_s.
\end{equation*}
Similarly, if $Z_t\equiv\theta(t,X_t)$ then $(1-V)Z=0$. Applying the above equality together with identity (\ref{eq:lemma:int of zcqv is zcqv:1}),
\begin{equation}\label{eq:pf:mufx for disjoint supports:1}\begin{split}
\int_0^t Z_{s-}\,dY_s
&=
\int_0^t Z_{s-}V_{s-}\,dY_s
=\sum_{s\le t}Z_{s-}V_{s-}Y_s\\
&=\sum_{s\le t}Z_{s-}Y_s
=\sum_{s\le t}J^X_s(\theta,f).
\end{split}\end{equation}
For the last equality, we have used the fact that the supports of $\theta$ and $f$ are disjoint in Definition \ref{def:J} to get $J^X_s(\theta,f)=Z_{s-}Y_s$.

Then, Lemma \ref{lemma:sum J is integrable} shows that $\int Z_-\,dY$ has integrable variation and, by the Doob-Meyer decomposition, can be written as $N+B$ for a \cadlag\ martingale $N$ and previsible process $B$ with integrable variation and $B_0=0$.
If $Y=M+A$ is decomposition (\ref{eqn:lemma:f is M plus A for nondif}) then $N-\int Z_-\,dM=\int Z_-\,dA-B$ is a previsible local martingale with zero continuous quadratic variation, and so is constant. Therefore, $\int Z_-\,dA=B$ has integrable variation. So $\theta^-$ is $\mu^X_f$-measurable and,
\begin{equation*}
\mu^X_f(\theta^-)=\E{B_\infty}=\E{\int_0^\infty Z_{s-}\,dY_s}=\E{\sum_{t>0}J^X_t(\theta,f)}=\tilde\mu^X_f(\theta).
\end{equation*}
The final equality just follows from Definition \ref{def:mutfx} of $\tilde\mu^X_f(\theta)$, noting that the $\mu_{[f,C]}(\theta)$ and $\mu_X$ terms are zero due to the supports of $\theta$ and $f$ being disjoint.
\end{proof}

We require the following quadratic covariations result from \citep{Lowther3}.
If $X$ is a semimartingale and $f\in\nicefunc$ then the quadratic variation of $f(t,X_t)$ satisfies
\begin{equation}\label{eqn:cov of f and g}
\left[f(\cdot,X_\cdot)\right]_t=\int_0^t f\pd{2}(s,X_{s})^2\,d[X]^c_s+\sum_{s\le t}(\Delta f(s,X_s))^2.
\end{equation}
See \citep{Lowther3} Theorem 1.4 and Lemma 1.5.
This allows us to prove the following integration by parts formula for $\tilde\mu^X_f$.
\begin{lemma}\label{lemma:int by parts for mutf}
Let $X$ be a class (DL) quasimartingale and $f,g\in\nicefuncK$. Then the quadratic covariation of $f(t,X_t)$ and $g(t,X_t)$ exists, has integrable variation and,
\begin{equation}\label{eqn:lemma:int by parts for mutf:1}
\tilde\mu^X_f(g) + \tilde\mu^X_g(f) + \E{\left[ f(\cdot,X_{\cdot}),g(\cdot,X_{\cdot})\right]_\infty}=0.
\end{equation}
\end{lemma}
\begin{proof}
Set $Y_t=f(t,X_t)$ and apply Definition \ref{def:mutfx},
\begin{equation}\label{eq:pf:int by parts for mutf:1}\begin{split}
\tilde\mu^X_f(f) ={}& \dblint(f\pd{2}C\pd{2}-f^-\pd{2}C^-\pd{2})\,d_tf\,dx-\dblint(f^-\pd{2})^2\,d_tC\,dx\\
&+\int\int_{-\infty}^xf^-\pd{2}(t,y)^2\,dy\,d\mu_X(t,x)\\
&+\E{\sum_{t>0}\int_{X_{t-}}^{X_t}((f-Y_t)f^-\pd{2}+(X_t-x)(f^-\pd{2})^2)\,dx}.
\end{split}\end{equation}
Substitute $(f^-\pd{2})^2$ for $\theta$ in equation (\ref{eqn:lemma:vol measure:1}) and subtract from (\ref{eq:pf:int by parts for mutf:1}),
\begin{equation}\label{eq:pf:int by parts for mutf:2}\begin{split}
\tilde\mu^X_f(f)
={}&
\dblint(f\pd{2}C\pd{2}-f^-\pd{2}C^-\pd{2})\,d_tf\,dx
-\frac{1}{2}\E{\int f\pd{2}(t,X_t)^2\,d[X]^c_t}\\
&+\E{\sum_{t>0}\int_{X_{t-}}^{X_t}(f-Y_t)f^-\pd{2}\,dx}
\end{split}\end{equation}
We require the following identity, which is simply an application of integration by parts. Here, the notation $\Delta_tf\equiv f-f^-$ is used.
\begin{equation}\label{eq:pf:int by parts for mutf:3}
\int_{X_{t-}}^{X_t}(f-Y_t)f^-\pd{2}\,dx+\frac{1}{2}\Delta Y_t^2
=\int\left(f\pd{2}1_{\{x<X_t\}}-f^-\pd{2}1_{\{x<X_{t-}\}}\right)\Delta_tf\,dx.
\end{equation}
Note that, as $f(t,x)$ is \cadlag\ in $t$, there can be only countably many times $t>0$ at which $f(t,x)\not= f^-(t,x)$. So, summing the right hand side of (\ref{eq:pf:int by parts for mutf:3}) over all $t>0$ gives
\begin{equation*}
\dblint\left(f\pd{2}1_{\{x<X_t\}}-f^-\pd{2}1_{\{x<X_{t-}\}}\right)d_tf\,dx.
\end{equation*}
As $f\pd{2}$ has bounded support, this integral is bounded. Then, sum (\ref{eq:pf:int by parts for mutf:3}) over $t>0$ and take expectations,
\begin{equation*}
\E{\sum_{t>0}\left(\int_{X_{t-}}^{X_t}(f-Y_t)f^-\pd{2}\,dx+\frac{1}{2}(\Delta Y_t)^2\right)}=
\dblint\left(f^-\pd{2}C^-\pd{2}-f\pd{2}C\pd{2}\right)d_tf\,dx.
\end{equation*}
Here, we have substituted in $-C\pd{2}$ and $-C^-\pd{2}$ for $\Prob{X_t>x}$ and $\Prob{X_{t-}>x}$ respectively. This equality can be added to equation (\ref{eq:pf:int by parts for mutf:2}) to get
\begin{equation*}
\tilde\mu^X_f(f)=-\frac{1}{2}\E{\int f\pd{2}(t,X_t)^2\,d[X]^c_t+\sum_{t>0}\Delta Y_t^2}=-\frac{1}{2}\E{[Y]_\infty}.
\end{equation*}
Here, equality (\ref{eqn:cov of f and g}) has been used. This proves the result in the case where $f=g$. The general result then follows from the linearity of (\ref{eqn:lemma:int by parts for mutf:1}) in both $f$ and $g$.
\end{proof}
This integration by parts formula can be used to exchange the roles of $\theta$ and $f$ in the expression  for $\mu^X_f(\theta^-)$.
\begin{lemma}\label{lemma:mufg mugf switch}
Let $X$ be a class (DL) quasimartingale, and $f,\theta\in\nicefuncK$. If $\theta\in\Dom(X)$ then $\theta^-$ is $\mu^X_f$-integrable. If, furthermore, $\mu^X_\theta(f^-)=\tilde\mu^X_\theta(f)$ then $\mu^X_f(\theta^-)=\tilde\mu^X_f(\theta)$.
\end{lemma}
\begin{proof}
Set $U_t=f(t,X_{t})$ and $V_t=\theta(t,X_{t})$ and use the decompositions $U=M+A$, $V=N+B$
where $M,N$ are local martingales and $A,B$ are previsible \zcqv\ processes. Integration by parts gives
\begin{equation}\label{eq:pf:mufg mugf switch:1}
\int_0^tU_{s-}\,dB_s+\int_0^tV_{s-}\,dA_s+[U,V]_t = U_tV_t-L_t
\end{equation}
where $L$ is the local martingale
\begin{equation*}
L_t=\int_0^t U_{s-}\,dN_s+\int_0^tV_{s-}\,dM_s.
\end{equation*}
However, by equation (\ref{eq:cov of X and A}), $[M,A]$ and $[N,B]$ are local martingales, so $\E{[U]}=\E{[M]+[A]}$ and $\E{[V]}=\E{[N]+[B]}$.
As Lemma \ref{lemma:int by parts for mutf} says that $[U]$ and $[V]$ are integrable, it follows that $[M]$ and $[N]$ are integrable, so $M$ and $N$ are square integrable martingales. Therefore $L$ is also a square integrable martingale. In particular, it is a martingale.
If $\theta$ is in $\Dom(X)$ then $\int U_-\,dB$ has integrable variation, and it follows from (\ref{eq:pf:mufg mugf switch:1}) that
\begin{equation*}
\sup_{t\ge 0}\left|\int_0^tV_{s-}\,dA_s\right|
\end{equation*}
is integrable, so $\theta^-$ is $\mu^X_f$-integrable. Take expectations of (\ref{eq:pf:mufg mugf switch:1}) and substitute in the definitions of $\mu^X_f$, $\mu^X_\theta$ to get
\begin{equation*}
\mu^X_\theta(f^-)+\mu^X_f(\theta^-)+\E{[U,V]_\infty}=0.
\end{equation*}
Combine this with equation (\ref{eqn:lemma:int by parts for mutf:1}) to eliminate the quadratic covariation term,
\begin{equation*}
\mu^X_\theta(f^-)+\mu^X_f(\theta^-)=\tilde\mu^X_\theta(f^-)+\tilde\mu^X_f(\theta)
\end{equation*}
and the result follows.
\end{proof}
Lemma \ref{lemma:mufg mugf switch} can be applied to calculate $\mu^X_f(\theta)$ for general $f\in\nicefunc$, in the case where $\theta$ is twice continuously differentiable.
\begin{lemma}\label{lemma:mux for theta smooth}
Let $X$ be a class (DL) quasimartingale and $f\in\nicefunc$. If $\theta\in\nicefuncK$ is twice continuously differentiable then $\theta$ is $\mu^X_f$-integrable and
$\mu^X_f(\theta)=\tilde\mu^X_f(\theta).$
\end{lemma}
\begin{proof}
By choosing any $g\in\nicefuncK$ such that $g=1$ on an open set containing the support of $\theta$, and writing $f=gf+(1-g)f$, then the linearity in $f$ means that we only need to prove the result separately for the cases where $f\in\nicefuncK$ and where $f$ and $\theta$ have disjoint support.

First, suppose that $f\in\nicefuncK$. Then Lemma \ref{lemma:muf is mutf for smooth f} says that $\theta\in\Dom(X)$ and $\mu^X_\theta(f^-)=\tilde\mu^X_\theta(f)$, and Lemma \ref{lemma:mufg mugf switch} gives the required result.
Finally, if $f$ and $\theta$ have disjoint support then the result follows from Lemma \ref{lemma:mufx for disjoint supports}.
\end{proof}

In order to prove Theorem \ref{thm:mux for gen proc} we just need to extend the previous result to non-differentiable $\theta$. The idea is simple enough -- just smooth $\theta$ by convolving with smooth functions and take limits of $\tilde\mu^X_f(\theta)$. However, it is rather tricky to show that these limits do in fact converge to the desired limits. We do this by smoothing $\theta(t,x)$ separately in $t$ and $x$.

\begin{proof}[Proof Of Theorem \ref{thm:mux for gen proc}]
First suppose that $\theta(t,x)$ is twice differentiable with respect to $x$ and $\theta\pd{2},\theta\pd{22}\in\nicefuncK$.
Then, choose any twice continuously differentiable $\alpha:(0,1)\rightarrow\reals$ with compact support and whose integral is equal to $1$. Define
\begin{equation*}
\theta^n(t,x)=\int_0^1\theta(t-s/n,x)\alpha(s)\,ds.
\end{equation*}
For $n$ large enough such that the support of $\theta$ is contained in $(1/n,\infty)\times\reals$, it follows that $\theta^n\in\nicefuncK$ is twice continuously differentiable, so Lemma \ref{lemma:mux for theta smooth} gives $\mu^X_f(\theta^n)=\tilde\mu^X_f(\theta^n)$. As $\theta^n\rightarrow\theta^-$, dominated convergence implies that $\mu^X_f(\theta^n)$ converges to $\mu^X_f(\theta^-)$.

We now show that $\tilde\mu^X_f(\theta^n)$ converges to $\tilde\mu^X_f(\theta)$.
As $(\theta^n)\pd{2}$ converges to $\theta^-\pd{2}$, dominated convergence implies that the drift and jump terms in Definition \ref{def:mutfx} do indeed converge. The remaining term is
\begin{equation*}
\mu_{[f,C]}(\theta^n)=\dblint f\pd{2}C\pd{2}\,d_t\theta^n\,dx -\dblint \theta^n\pd{2} f\pd{2}^-\,d_tC\,dx-\dblint \theta^n\pd{2} C\pd{2}^-\,d_tf\,dx.
\end{equation*}
Again by dominated convergence, the last two terms on the right hand side will converge. We look at the first term, which can be rearranged to get
\begin{equation}\label{eq:pf:mux for gen proc:1}
\dblint f\pd{2}C\pd{2}\,d_t\theta^n\,dx=\dblint g\,d_t\theta^n\pd{2}\,dx=\dblint\int_{0}^1 g(t+s/n,x)\,ds\,d_t\theta\pd{2}\,dx
\end{equation}
where
\begin{equation}\label{eq:pf:mux for gen proc:2}\begin{split}
g(t,x)&\equiv \int_x^{\infty}f\pd{2}(t,y)C\pd{2}(t,y)\,dy=-\int_x^\infty f\pd{2}(t,y)\Prob{X_t>y}\,dy\\
&=\E{1_{\{X_t>x\}}(f(t,x)-f(t,X_t))}.
\end{split}\end{equation}
As $f$ is Lipschitz continuous in $x$ and $X_t$ is right-continuous and uniformly integrable on bounded intervals, it follows that $g(t,x)$ is right-continuous in $t$. So, taking limits in (\ref{eq:pf:mux for gen proc:1}) and applying dominated convergence,
\begin{equation*}
\dblint f\pd{2}C\pd{2}\,d_t\theta^n\,dx\rightarrow\dblint g\,d_t\theta\pd{2}\,dx=\dblint f\pd{2}C\pd{2}\,d_t\theta\,dx
\end{equation*}
as $n\rightarrow\infty$.
Therefore, $\tilde\mu^X_f(\theta^n)=\mu^X_f(\theta^n)$ converges to $\tilde\mu^X_f(\theta)$, and it follows that $\mu^X_f(\theta^-)=\tilde\mu^X_f(\theta)$.

Now choose any $\theta\in\nicefuncK$ and twice continuously differentiable $\alpha:\reals\rightarrow\reals$ with compact support and which integrates to $1$. Set
\begin{equation}\label{eq:pf:mux for gen proc:3}
\theta^n(t,x)=\int\theta(t,x+y/n)\alpha(y)\,dy.
\end{equation}
Then, $\theta^n(t,x)$ is twice differentiable in $x$ and $\theta^n$, $\theta^n\pd{2}$ and $\theta^n\pd{22}$ are all in $\nicefuncK$. The argument above shows that $\mu^X_f((\theta^n)^-)=\tilde\mu^X_f(\theta^n)$.
Choose any $(t,x)\in\halfplane$ such that $f(t,x)$ and $C(t,x)$ are differentiable in $x$. Then, $\Prob{X_t=x}=0$ and if $g$ is defined by equation (\ref{eq:pf:mux for gen proc:2}) it follows that $g(t,x)$ is differentiable with respect to $x$ and,
\begin{equation*}
g\pd{2}(t,x)=\E{1_{\{X_t>x\}}f\pd{2}(t,x)}=-f\pd{2}(t,x)C\pd{2}(t,x).
\end{equation*}
By Lemma \ref{lemma:nicefunc ae diff}, $f(t,x)$ and $C(t,x)$ are differentiable in $x$ almost everywhere with respect to the measure $\dblint\cdot\,|d_t\theta|\,dx$, and dominated convergence gives
\begin{equation*}
\dblint f\pd{2}C\pd{2}\,d_t\theta^n\,dx=-\dblint \int g\pd{2}(t,x-y/n)\,dy\,d_t\theta\,dx\rightarrow\dblint f\pd{2}C\pd{2}\,d_t\theta\,dx
\end{equation*}
as $n\rightarrow\infty$.
Also, from the definition (\ref{eq:pf:mux for gen proc:3}) above for $\theta^n$ it follows that $\theta^n\pd{2}(t-,x)$ converges to $\theta\pd{2}(t-,x)$ as $n\rightarrow\infty$ whenever $\theta^-(t,x)$ is differentiable with respect to $x$.
So, as above, dominated convergence implies that $\mu^X_f((\theta^n)^-)$ tends to $\mu^X_f(\theta^-)$ and $\tilde\mu^X_f(\theta^n)$ tends to $\tilde\mu^X_f(\theta)$. So, $\tilde\mu^X_f(\theta)=\mu^X_f(\theta^-)$ as required.
\end{proof}

\section{The Time Reversed Conditional Variation}
\label{sec:reversed condvar}

In the previous section the expression for the drift of $f(t,X_t)$ given by Theorem \ref{thm:mux for gen proc} was proved, under the condition that $f$ is in $\Dom(X)$.
However, it still remains to show that the conditions specified in Theorem \ref{thm:mufc defines drift of acd process} are sufficient for $f$ to be in $\Dom(X)$. This direction is complicated by the problem that $f(t,X_t)$ cannot be assumed a priori to be a semimartingale.
The idea is to bound the conditional variation, and show that $f(t,X_t)$ is locally a quasimartingale. Then, Rao's decomposition can be used to decompose it into local martingale and finite variation terms.
Unfortunately, the conditions on $f$ are too weak for this to be an easy route, and instead we will bound the conditional variation of the drift term $A$ from decomposition (\ref{eqn:lemma:f is M plus A for nondif}) under time reversal. This will be enough to show that $A$ is a finite variation process and, consequently, $f(t,X_t)$ will be a semimartingale.
The following time reversed natural filtration of $X$ will be used.
\begin{equation*}
\setsG^X_t=\sigma\left(\{X_s:s\ge t\}\cup\{A\in\setsF:\Prob{A}=0\}\right).
\end{equation*}
Note that for times $s\le t$, then $\setsG^X_s\supseteq\setsG^X_t$. The conditional variation with respect to this filtration is defined as follows.
\begin{definition}\label{def:rev cond var}
Let $X,A$ be real valued processes such that $\E{|A_t|}<\infty$ for every $t\ge 0$. Define
\begin{equation*}
V^r_X(A)=\sup\E{\sum_{k=1}^nZ_k(A_{t_k}-A_{t_{k-1}})},
\end{equation*}
where the supremum is taken over all sequences $0\le t_0\le t_1\le \cdots\le t_n$ and all $\setsG^X_{t_{k}}$-measurable random variables $|Z_k|\le 1$.
\end{definition}
In the current section we show that bounding the time reversed conditional variation $V^r_X(A)$ is enough to bound the variation of $A$. Then, in the following section, this will be applied to the proofs of Theorems \ref{thm:abs drift is abs measure} and \ref{thm:mufc defines drift of acd process}.
The main result used is the following.
\begin{lemma}\label{lemma:bdd rev cond var gives int var}
Let $X$ be a \cadlag\ real valued process and $A$ be a \zcqv\ process such that $\sup_{t\ge 0}|A_t|$ is integrable and $V^r_X(A)<\infty$.

Suppose furthermore that $A_t-A_s$ is $\setsG^X_s$-measurable for all $t>s$ and that there exists a measurable $u:\halfplane\rightarrow\reals$ such that $\Delta A_t=u(t,X_t)$ for all $t> 0$.
Then, $A$ has integrable variation satisfying
\begin{equation*}
\E{\int\,|dA|}\le V^r_X(A).
\end{equation*}
\end{lemma}
\begin{proof}
Pick any $T>0$ and let $\setsF^r_\cdot$ be the filtration and $B$ be the \cadlag\ process defined as
\begin{equation*}
\setsF^r_t = \setsG^X_{(T-t)-},\ 
B_t=A_{T-}-A_{(T-t)-}.
\end{equation*}
In order that this makes sense for $t\ge T$, we set $\setsF_t=\setsF_0$ and $A_t=A_0$ for $t<0$.
The condition that $A_t-A_s$ is $\setsG^X_s$-measurable for all $t>s$ ensures that $B$ is $\setsF^r_\cdot$-adapted.
By the condition of the lemma, $\Delta B_t=\Delta A_{T-t}=u(T-t,X_{T-t})$ is a measurable function of the left-continuous and $\setsF^r_\cdot$-adapted process $X_{T-t}$. So, $\Delta B$ and hence $B$ are $\setsF^r_\cdot$-previsible.

Now choose any $0\le t_0\le t_1\le\cdots\le t_n$ and random variables $|Z_k|\le 1$ such that $Z_k$ is $\setsF^r_{t_{k-1}}$-measurable.
Setting $s^\alpha_k=T-t_{n-k}-\alpha$ for any real $\alpha>0$,
\begin{equation*}
B_{t_k}-B_{t_{k-1}}=\lim_{\alpha\rightarrow 0}(A_{s^\alpha_{n+1-k}}-A_{s^\alpha_{n-k}}).
\end{equation*}
Also, $\setsF^r_{t_{k-1}}\subseteq\setsG_{s^\alpha_{n+1-k}}$, so $Z_k$ is $\setsG_{s^\alpha_{n+1-k}}$-measurable, and dominated convergence gives,
\begin{equation*}
\sum_{k=1}^n\E{Z_k(B_{t_k}-B_{t_{k-1}})}=\lim_{\alpha\rightarrow 0}\sum_{k=1}^n\E{Z_{n+1-k}(A_{s^\alpha_k}-A_{s^\alpha_{k-1}})}
\le V^r_X(A).
\end{equation*}
So $B$ is an $\setsF^r_\cdot$-quasimartingale with conditional variation bounded by $V^r_X(A)$, and Rao's decomposition (Lemma \ref{lemma:Rao decomposition}) can be applied,
\begin{equation}\label{eqn:pf:mufc defines drift of acd process:2}
B_t=N_t+C_t.
\end{equation}
Here, $N$ is an $\setsF^r_\cdot$-local martingale and $C$ is a \cadlag\ $\setsF^r_\cdot$-previsible process with
\begin{equation*}
\E{\int_0^T\,|dC|}\le V^r_X(A).
\end{equation*}
However, $B$ and $C$, and therefore $N$, are $\setsF^r_\cdot$-previsible processes, so $N$ must be continuous.
Also, as $A$ is a \zcqv\ process, it follows that $N+C$ and therefore $N$ also have zero continuous quadratic variation. Furthermore, the continuity of $N$ implies that $[N]=0$ and, as $N$ is a martingale it must be constant. So,
\begin{equation*}
A_t=A_{T-}-B_{(T-t)-}=A_{T-}-N_0-C_{(T-t)-}
\end{equation*}
for $t<T$.
Finally,
\begin{equation*}
\E{\int_0^{T-}\,|dA|}=\E{\int_0^T\,|dC|}\le V^r_X(A).
\end{equation*}
The result follows by letting $T$ increase to infinity.
\end{proof}
The remainder of this section will be used to show that the conditions required for Lemma \ref{lemma:bdd rev cond var gives int var} to be applied are satisfied in the cases we need, starting with the following result which expresses $\Delta A$ as a function of $X$.
\begin{lemma}\label{lemma:DA is func of X}
Let $X$ be an adapted and quasi-left-continuous process and $f\in\nicefunc$. Suppose that decomposition (\ref{eqn:lemma:f is M plus A for nondif}) exists, and is $f(t,X_t)=M_t+A_t$. Then,
\begin{equation*}
\Delta A_t = f(t,X_t)-f(t-,X_t)
\end{equation*}
for every $t>0$.
\end{lemma}
\begin{proof}
First, choose any previsible stopping time $\tau>0$. The quasi-left-continuity of $X$ gives $X_{\tau}=X_{\tau-}$. As $A$ is previsible, $\E{\Delta A_{\tau}|\setsF_{\tau-}}=\Delta A_\tau$ and, as $M$ is a local martingale, $\E{\Delta M_\tau|\setsF_{\tau-}}=0$. So,
\begin{equation*}
\Delta A_\tau=\E{f(\tau,X_{\tau})-f(\tau-,X_{\tau-})|\setsF_{\tau-}}=f(\tau,X_{\tau-})-f(\tau-,X_{\tau-}).
\end{equation*}
As $X$ is adapted, both sides of this equality are previsible. So, by previsible section, it remains true for all times $\tau>0$.
Now, as $f(t,x)$ is jointly continuous in $x$ and \cadlag\ in $t$, there exists a countable $S\subset\reals$ such that $f(t,x)=f^-(t,x)$ for all $t\not\in S$. Also, the quasi-left-continuity of $X$ implies that it is continuous in probability, so $\Delta X_t=0$ for every $t\in S$.
Then,
\begin{equation*}\begin{split}
\Delta A_t &= 1_{\{t\in S\}}\left( f(t,X_{t-})-f(t-,X_{t-})\right)=1_{\{t\in S\}}\left( f(t,X_{t})-f(t-,X_{t})\right)\\
&=f(t,X_{t})-f(t-,X_{t})
\end{split}\end{equation*}
as required.
\end{proof}
Another condition required by Lemma \ref{lemma:bdd rev cond var gives int var} is that $A_t-A_s$ should be $\setsG^X_s$-measurable whenever $t>s$. This will follow by applying the following lemma to decomposition (\ref{eqn:lemma:f is M plus A for nondif}).
\begin{lemma}\label{lemma:decomp is adapted}
Let $X$ be a Markov process and $Y=M+A$ where $M$ is a local martingale and $A$ is a previsible \zcqv\ process with $A_0=0$. Then, $Y$ is adapted to the natural filtration, $\setsF^X_\cdot$, of $X$ if and only if $A$ and $M$ are $\setsF^X_\cdot$-adapted.
\end{lemma}
\begin{proof}
If $M$ and $A$ are $\setsF^X_\cdot$-adapted then clearly $Y$ is also $\setsF^X_\cdot$-adapted. Conversely, suppose that $Y$ is $\setsF^X$-adapted.
For the moment, let us suppose that $[Y]_\infty$ is integrable.
As identity (\ref{eq:cov of X and A}) says that $[M,A]$ is a local martingale, monotone convergence gives $\E{[M]_\infty+[A]_\infty}=\E{[Y]_\infty}<\infty.$ So, $M$ is a uniformly square integrable martingale.

Now let $N$ be the \cadlag\ process $N_t=\E{M_t|\setsF^X_\infty}$. As $M$ is adapted, the Markov property for $X$ gives $N_t=\E{M_t|\setsF^X_t}$, so $N$ is a uniformly square integrable $\setsF^X_\cdot$-martingale. Set $\tilde Y=Y-N$ and $\tilde M=M-N$.

If, for every $n\in\nat$, we let $P_n$ be the partition of $\reals$ given by $\tau^{P_n}_k=k/n$ then, letting $[\tilde M]^{P_n}=[\tilde M,\tilde M]^{P_n}$ be the approximation to the quadratic variation of $\tilde M$ given by (\ref{eqn:qv approx def}),
\begin{equation*}
[\tilde M]^{P_n}_t-[\tilde M]_t=2\int_0^t\sum_{k=1}^\infty 1_{\{\tau^{P_n}_{k-1}< s\le\tau^{P_n}_k\}}(\tilde M_{s-}-\tilde M_{\tau^{P_n}_{k-1}})\,d\tilde M_s
\end{equation*}
is a local martingale. So, by monotone convergence, the expected values of the sequence $[\tilde M]^{P_n}_t$ are bounded by the expected value of $[\tilde M]_t$. Therefore, $([\tilde M]^{P_n}_t)^{1/2}$ is uniformly integrable over all $n\in\nat$.
As the Cauchy-Schwarz inequality implies that $|[\tilde Y,\tilde M]^{P_n}_t|$ is bounded by $([\tilde Y]^{P_n}_t[\tilde M]^{P_n}_t)^{1/2}$, we can take limits
\begin{equation*}\begin{split}
&\E{[\tilde Y,\tilde M]_t|\setsF^X_\infty}=\lim_{n\rightarrow\infty}\E{[\tilde Y,\tilde M]_t^{P_n}|\setsF^X_\infty}\\
&=\lim_{n\rightarrow\infty}\sum_{k=1}^\infty (\tilde Y_{t^{P_n}_k\wedge t}-\tilde Y_{t^{P_n}_{k-1}\wedge t})\E{\tilde M_{t^{P_n}_k\wedge t}-\tilde M_{t^{P_n}_{k-1}\wedge t}|\setsF^X_\infty}=0.
\end{split}\end{equation*}
Furthermore, $[A,\tilde M]$ is a local martingale with integrable variation bounded by $([A]_\infty[\tilde M]_\infty)^{1/2}$, so it is a proper martingale. Therefore,
\begin{equation*}
\E{[\tilde M]_t}=\E{[\tilde Y,\tilde M]_t}-\E{[A,\tilde M]_t}=0,
\end{equation*}
so, $M=N$ and $A=Y-N$ are $\setsF^X_\cdot$-adapted.

Now, returning to the general case, we can define stopping times
\begin{equation*}
\tau_m=\inf\left\{t\in\reals_+:[Y]_t\ge m\right\}.
\end{equation*}
Letting $Y^{\tau_m-}_t$ be the pre-stopped process equal to $Y_t$ for $t<\tau_m$ and $Y_{\tau_m-}$ for $t\ge\tau_m$, it follows that $[Y^{\tau_m-}]_\infty<m$.

As $Y$ is a sum of a local martingale and a \cadlag\ previsible (hence locally bounded) process, $Y^*_t\equiv\sup_{s\le t}|Y_s|$ will be locally integrable.
So, $Z^m\equiv1_{\slbrack\tau_m,\infty\slbrack}\Delta Y_{\tau_m}$ has locally integrable variation, and the Doob-Meyer decomposition gives $Z^m=N^m+B^m$ for a local martingale $N$ and \cadlag\ previsible process $B^m$ with locally integrable variation, and we can take $N^m_0=B^m_0=0$. Then, the argument above can be applied to $Y^{\tau_m-}=(M^{\tau_m}-N^m)+(A^{\tau_m}-B^m)$ to see that $A^{\tau_m}-B^m$ is $\setsF^X_\cdot$-adapted.

Finally, choosing any $t\ge 0$ and stopping time $\sigma$ such that $Y^*_\sigma$ is integrable,
\begin{equation*}\begin{split}
\E{1_{\{t<\sigma\}}|B^m_t|} &\le \E{\int_0^{t\wedge \sigma}\,|dB^m|}\le \E{1_{\{t\wedge\sigma\ge\tau_m \}}|\Delta Y_{\tau_m}|}\\ &\le 2\E{1_{\{t\ge\tau_m \}}Y^*_\sigma}\rightarrow 0
\end{split}\end{equation*}
as $m\rightarrow\infty$.
As $Y^*$ is locally integrable, $\Prob{\sigma<t}$ can be made as small as we like, showing that $B^m_t\rightarrow 0$ in probability as $m\rightarrow\infty$. So $A_t=\lim_{m\rightarrow\infty}(A^{\tau_m}_t-B^m_t)$ is $\setsF^X_\cdot$-adapted.
\end{proof}
Finally for this section, we apply the previous result to provide a sufficient condition for $A_t-A_s$ to be $\setsG^X_s$-measurable, as required by Lemma \ref{lemma:bdd rev cond var gives int var}.
\begin{corollary}\label{cor: diff A is G meas}
Let $X$ be a Markov process and $Y=M+A$ where $M$ is a local martingale and $A$ is a \cadlag\ previsible process with zero continuous quadratic variation.
If $Y_t-Y_s$ is $\setsG^X_s$-measurable for every $t>s$ then the same holds for $A$.
\end{corollary}
\begin{proof}
Fix any $s\ge 0$, and set
\begin{align*}
&\tilde X_t=1_{\{t\ge s\}}X_t,\ \tilde Y=1_{\{t\ge s\}}(Y_t-Y_s),\\
&\tilde M_t=1_{\{t\ge s\}}(M_t-M_s),\ \tilde A_t=1_{\{t\ge s\}}(A_t-A_s).
\end{align*}
Applying Lemma \ref{lemma:decomp is adapted} to $\tilde Y=\tilde M+\tilde A$, we see that $\tilde A$ is $\setsF^{\tilde X}_\cdot$-adapted. In particular, for any $t\ge s$, $A_t-A_s=\tilde A_t$ is measurable with respect to $\setsF^{\tilde X}_t\subseteq\setsG^X_s$.
\end{proof}

\section{Proof of Theorems \ref{thm:abs drift is abs measure} and \ref{thm:mufc defines drift of acd process}}
\label{sec:pf of 2.5 and 2.6}

We now apply Lemma  \ref{lemma:bdd rev cond var gives int var} to the proof of Theorem \ref{thm:abs drift is abs measure}.

\begin{proof}[ Proof of Theorem \ref{thm:abs drift is abs measure}]
Choose any bounded and measurable $\theta:\halfplane\rightarrow\reals$ with bounded support and set
\begin{equation*}
B_t=\int_0^t \theta(s,X_{s-})\,dA_s.
\end{equation*}
We will show that the conditions required to apply Lemma \ref{lemma:bdd rev cond var gives int var} to $B$ are satisfied. First, as $f\in\Dom(X)$, $B$ has integrable variation. In particular, it has zero continuous quadratic variation and is integrable.
Lemma \ref{cor: diff A is G meas} shows that $A_t-A_s$ is $\setsG^X_s$-measurable for every $t>s$, so $B_t-B_s$ will also be $\setsG^X_s$-measurable.
The quasi-left-continuity of $X$ implies that $X_t=X_{t-}$ whenever $\Delta A\not=0$. So, applying Lemma \ref{lemma:DA is func of X},
\begin{equation*}
\Delta B_t=\theta(t,X_t)\Delta A_t=\theta(t,X_t)(f(t,X_t)-f(t-,X_t)).
\end{equation*}
We now bound $V^r_X(B)$. Choose any $T>0$ and measurable $g:\reals\rightarrow[-1,1]$. Then, as $X$ is strong Markov there exists a measurable $h:\halfplane\rightarrow[-1,1]$ such that
\begin{equation*}
h(\tau,X_\tau)=\E{g(X_T)|\setsF_\tau}
\end{equation*}
for every stopping time $\tau\le T$. For every $t\le T$, optional projection and the quasi-left-continuity of $X$ gives,
\begin{equation*}\begin{split}
\E{g(X_T)(B_T-B_t)}
&=\E{\int_t^T h(s,X_{s-})\theta(s,X_{s-})\,dA_s}\\
&=\mu^X_f(1_{(t,T]}h\theta)\le |\mu^X_f|(1_{(t,T]}|\theta|).
\end{split}\end{equation*}
So, Definition \ref{def:rev cond var} gives
$V^r_X(B)\le|\mu^X_f|(|\theta|)$
and, finally, Lemma \ref{lemma:bdd rev cond var gives int var} can be applied,
\begin{equation*}
\E{\int|\theta(t,X_{t-})|\,|dA_t|}=\E{\int\,|dB|}\le |\mu^X_f|(|\theta|)
\end{equation*}
as required.
\end{proof}

The proof of Theorem \ref{thm:mufc defines drift of acd process} will require the following result, generalizing the expression for $\mu^X_f(\theta)$ to non-smooth $\theta$ and all $f$ in $\nicefunc$. The proof follows in the same way as for Lemma \ref{lemma:mux for theta smooth}.

\begin{lemma}\label{lemma:mux for theta nonsmooth}
Let $X$ be a class (DL) quasimartingale and $f\in\nicefunc$. If $\theta\in\nicefuncK$ and $\theta\in\Dom(X)$, then $\theta^-$ is $\mu^X_f$-integrable and
$\mu^X_f(\theta^-)=\tilde\mu^X_f(\theta).$
\end{lemma}
\begin{proof}
As in the proof of Lemma \ref{lemma:mux for theta smooth}, the linearity in $f$ means that it is enough to prove the result separately for the cases where $f\in\nicefuncK$ and where $f$ and $\theta$ have disjoint support.
First, if $f\in\nicefuncK$, Theorem \ref{thm:mux for gen proc} gives $\mu^X_{\theta}(f^-)=\tilde\mu^X_{\theta}(f)$. Then, the result follows from Lemma \ref{lemma:mufg mugf switch}.
Finally, if $f$ and $\theta$ have disjoint support, the result is given by Lemma \ref{lemma:mufx for disjoint supports}.
\end{proof}

We also need to know when the product of functions $f,g$ in $\Dom(X)$ is itself in $\Dom(X)$.
This will be true whenever $fg$ is Lipschitz continuous, although we only require the case where either $f$ or $g$ has compact support. 

\begin{lemma}\label{lemma:products in DomX}
Let $X$ be a class (DL) quasimartingale and $f,g\in\Dom(X)$. If either $f$ or $g$ is in $\nicefuncK$ then $fg\in\Dom(X)$.
\end{lemma}
\begin{proof}
Without loss of generality, suppose that $g\in\nicefuncK$. As with the proofs of Lemmas \ref{lemma:mux for theta smooth} and \ref{lemma:mux for theta nonsmooth}, the linearity in $f$ means that it is enough to prove the result separately for the cases where $f,g$ have disjoint supports and where $f\in\nicefuncK$.

If they have disjoint support then $fg=0\in\Dom(X)$. If $f\in\nicefuncK$ then, setting $f_t\equiv f(t,X_t)$ and $g_t\equiv g(t,X_t)$, Lemma \ref{lemma:mufg mugf switch} says that the quadratic covariation $[f,g]$ has integrable variation. If $f_t=M_t+A_t$ and $g_t=\tilde M_t+\tilde A_t$ are decomposition (\ref{eqn:decomposition of f in domain of X}) then integration by parts,
\begin{equation*}\begin{split}
f_tg_t&=\int_0^t f_{s-}\,dg_s +\int_0^t g_{s-}\,df_s +[f,g]_t\\
&=\int_0^tf_{s-}\,d\tilde M_s+\int_0^tg_{s-}dM_s + \int_0^tf_{s-}\,d\tilde A_s+\int_0^tg_{s-}\,dA_s+[f,g]_t,
\end{split}\end{equation*}
shows that $f_tg_t$ decomposes into local martingale and integrable variation terms, so $fg\in\Dom(X)$.
\end{proof}

We finish with the proof of Theorem \ref{thm:mufc defines drift of acd process}, which uses Lemma \ref{lemma:bdd rev cond var gives int var} together with integrals of \zcqv\ processes studied in Section \ref{sec:gen drift}.

Lemma \ref{thm:f is convex dec} will also be required. If $X$ is a continuous and strong Markov martingale, $T>0$ and $g:\reals\rightarrow\reals$ is convex and Lipschitz continuous then this says that there exists an $h\in\nicefunc$ such that
\begin{equation}\label{eqn:h is exp g 3}
h(t,X_t)=\E{g(X_T)|\setsF_t}
\end{equation}
for every $t<T$. By linearity, this remains true when $g$ is a difference of convex and Lipschitz continuous functions. If, furthermore, $|g|\le 1$ then $h$ can be replaced by $(h\vee -1)\wedge 1$, which will also be in $\nicefunc$.

\begin{proof}[Proof Of Theorem \ref{thm:mufc defines drift of acd process}]
First, if $f\in\Dom(X)$ then the result follows from Theorem \ref{thm:mux for gen proc}. Conversely, suppose that there is a local signed measure $\mu$ satisfying $\mu_{[f,C]}(\theta)=\mu(\theta^-)$ for all $\theta\in\nicefuncK$. As $X$ is continuous and a martingale, this says that $\tilde\mu^X_f(\theta)=\mu(\theta^-)$.
Let $f(t,X_t)=M_t+A_t$ be decomposition (\ref{eqn:lemma:f is M plus A for nondif}) and choose any twice continuously differentiable $\theta\in\nicefuncK$. Set
\begin{equation*}
B_t=\int_0^t \theta(s,X_{s})\,dA_s.
\end{equation*}
We shall show that the necessary conditions for $B$ are satisfied in order for Lemma \ref{lemma:bdd rev cond var gives int var} to be applied.
First, by Lemma \ref{lemma:int of zcqv is zcqv}, $B$ is a previsible \zcqv\ process.
Also, Lemma \ref{cor: diff A is G meas} says that $A_t-A_s$ is $\setsG^X_s$-measurable for every $t>s$, and it follows that $B_t-B_s$ is $\setsG^X_s$-measurable. Furthermore, Lemma \ref{lemma:DA is func of X} gives
\begin{equation*}
\Delta B_t=\theta(t,X_t)\Delta A_t=\theta(t,X_t)(f(t,X_t)-f^-(t,X_t)).
\end{equation*}
It just remains to bound $V^r_X(B)$. Pick any $T>0$ let $g:\reals\rightarrow[-1,1]$ be a difference of Lipschitz continuous convex functions. Then, there exists an $h\in\nicefunc$ such that $|h|\le 1$ and (\ref{eqn:h is exp g 3}) is satisfied. Setting $N_t=h(t,X_t)$,
\begin{equation}\label{eq:pf:mufc for cts mgales:2}
h(t,X_t)(B_t-B_s) = \int_s^t h^-(u,X_u)\theta(u,X_u)\,dA_u + \int_s^t B_u\,dN_u.
\end{equation}
for all $t\ge s$ in the interval $[0,T)$. Here, identity (\ref{eq:lemma:int of zcqv is zcqv:1}) and the integration by parts formula (\ref{eqn:int wrt prev zcqv}) have been applied. The final term on the right hand side is a local martingale (with parameter $t\in[s,T)$).
As Lemma \ref{lemma:products in DomX} says that $h\theta\in\Dom(X)$, we can take expectations of (\ref{eq:pf:mufc for cts mgales:2}) and apply Lemma \ref{lemma:mux for theta nonsmooth} to get
\begin{equation*}
\E{g(X_T)(B_t-B_s)}=\tilde\mu^X_f(1_{[s,t)}h\theta)=\mu(1_{(s,t]}h^-\theta)\le |\mu|(1_{(s,t]}|\theta|).
\end{equation*}
Letting $T$ decrease to $t$ then gives
\begin{equation*}
\E{g(X_t)(B_t-B_s)}\le|\mu|(1_{(s,t]}|\theta|)
\end{equation*}
for all $t>s>0$. By the monotone class lemma, this extends to all measurable $g:\reals\rightarrow[-1,1]$, so for any $\setsG^X_t$-measurable random variable $|Z|\le 1$, the Markov property for $X$ implies that
\begin{equation*}
\E{Z(B_t-B_s)}=\E{\E{Z|X_t}(B_t-B_s)}\le |\mu|(1_{(s,t]}|\theta|).
\end{equation*}
Putting this into Definition \ref{def:rev cond var} of $V^r_X(B)$, shows that
$V^r_X(B)\le|\mu|(|\theta|).$
Therefore Lemma \ref{lemma:bdd rev cond var gives int var} can be applied, so $B$ has integrable variation, with expectation bounded by $|\mu|(|\theta|)$.
So,
\begin{equation*}
\E{\int|\theta(s,X_{s-})|\,|dA_s|}\le |\mu|(|\theta|)<\infty.
\end{equation*}
By monotone convergence, this extends to all continuous $\theta$ with bounded support. So $A$ satisfies the requirements of decomposition (\ref{eqn:decomposition of f in domain of X}), and $f$ is in $\Dom(X)$.
\end{proof}

\bibliography{backward.bbl}
\bibliographystyle{plain}

\end{document}